\documentclass[11pt]{amsart}
\usepackage{amsbsy,amsfonts}
\usepackage{latexsym,eucal,amsmath,amsthm}
\usepackage{amssymb}
\newcommand{\db}{\mathbb }
\newcommand{\dbC}{{\db C}}
\newcommand{\dbR}{{\db R}}

\newcommand{\dbZ}{{\db Z}}
\newcommand{\dbN}{{\db N}}

\newcommand{\F}{\mbox{{\em F}}}

\newcommand{\x}{\xi}
\newcommand{\z}{\zeta}
\newcommand{\la}{\lambda}

\def\Re{\mbox{\rm Re}}

\theoremstyle{plain}
   \newtheorem{theorem}[subsection]{Theorem}
   
   \newtheorem{proposition}[subsection]{Proposition}
   \newtheorem{lemma}[subsection]{Lemma}

\theoremstyle{remark}
   \newtheorem{remark}[subsection]{Remark}

\theoremstyle{definition}
   \newtheorem{definition}[subsection]{Definition}

\begin{document}
\date{2 November 2001}
\title{On Solutions for the Kadomtsev-Petviashvili I Equation}
\author{J. Colliander}
\thanks{J.E.C. was supported in part by  N.S.F. Grant DMS 0100595.}
\address{\small University of Toronto}

\author{C. Kenig}
\thanks{C.K. was supported in part by N.S.F Grant DMS 9500725}
\address{\small University of Chicago}

\author{G. Staffilani}
\thanks{G.S. was supported in part by N.S.F. Grant DMS 9800879,
the Terman Award and a grant by the Sloan Foundation.}
\address{\small Brown University}

\begin{abstract}
Oscillatory integral techniques are used to
study the well-posedness of the KP-I equation for
initial data that are small with respect to the norm of a
weighted Sobolev
space involving derivatives of total order no larger than $2$.
\end{abstract}

\maketitle

\label{intro}\section{Introduction}
We consider the initial value problem (IVP) for the
Kadomtsev-Petviashvili equation
\begin{equation}
\left\{ \begin{array}{l}
\partial_x(\partial_t u + \partial_x^3 u
+\beta\partial_xu^2)+ \gamma\partial_y^2u= 0, \\
u(x,0) = u_0(x), \hspace{1.5cm}(x,y) \in \dbR^2, \,  t \in \dbR,
\end{array}\right.
\label{ivp1}\end{equation}
where $u=u(t,x,y)$ is a scalar unknown function, $\beta$  is a real
constant and  $\gamma=\pm 1$.  The KP equation
models \cite{KP} 
the propagation along the $x$-axis of nonlinear dispersive long
waves on the surface of a fluid with a slow variation along the
$y$-axis. KP arises as a universal model in wave propagation and may
be viewed as a 2d generalization of the KdV equation. 
If we denote with $\partial_{x}^{-1}$ the antiderivative with
respect to the variable $x$, then we can rewrite the
evolution equation in (\ref{ivp1}) as
\begin{equation}
\partial_t u + \partial_x^3 u + \gamma\partial_x^{-1}\partial_y^2u
+\beta\partial_xu^2= 0.
\label{equation}\end{equation}
We will be using \eqref{equation} in the rest of the paper. This
equation is of {\em dispersive} type and the
strength of the dispersive effect depends on the sign of $\gamma$. To
see this we recall that the solution of the linear problem associated
to (\ref{equation}) can be written as the oscillatory integral
\begin{equation}
   U(t)u_{0}(x,y)=\int_{\dbR^{2}}e^{i(t\phi(\x,\lambda)+x\xi+y\lambda)}
\widehat{u_{0}}(\xi,\lambda)\, d\x\, d\lambda
\label{lsolution}\end{equation}
where $\phi(\x,\lambda)=\x^{3}-\gamma
\frac{\lambda^{2}}{\xi}$ is often called the {\em dispersive function}
or {\em dispersion relation}.
Hence one can
interpret the solution $U(t)u_{0}$ as the adjoint of the restriction
of the Fourier transform on the surface
$S=\{(\lambda,\xi)/ \lambda=\phi(\xi)\}$. It is well known, in
particular from the work of Stein \cite{Stein}, that the curvature
of the surface plays an important role in obtaining good
estimates for the restricted Fourier transform. The dispersive
function defines also the intensity of the smoothing effect. 
Kenig, Ponce and Vega \cite{KPV1} proved that if
$$u(x_{1},\ldots.,x_{n},t)=\int_{\dbR^{n}}
e^{i(t\phi(\x_{1},\ldots.,\x_{n})+\sum_{i=1}^{n}x_{i}\x_{i})}
\widehat{u_{0}} (\xi_1 , \dots, \xi_n) \, d\x_{1}\ldots.\, d\x_{n}$$
for a generic
dispersive function $\phi(\x_{1},\ldots, \x_{n})$ satisfying
$|\nabla\phi(\xi)|\geq C \sum_{i=1}^{n}
|\xi_{i}|^{\delta_{i}}$, for $\delta_{i}>0$, it follows that if
$u_{0} \in L^{2}$ then
$\partial_{x_{i}}^{\delta_{i}/2}u$ is a function in some $L^{p}$
space. If we go back to (\ref{equation}) we see that
\begin{equation}
|\nabla_{(\xi,\nu)}\phi(\xi,\lambda)|\geq C \left\{ \begin{array}{ll}
|\xi| &\mbox{ if } \gamma=-1\\
|\xi|^{2} &\mbox{ if } \gamma =1.
\end{array}\right.
\label{gradient}\end{equation}
We call the equation (\ref{equation}) KP-I if $\gamma=-1$ and KP-II
if  $\gamma=1$. It is clear then from (\ref{gradient})
that the linear solution for the KP-I equation gains in general no
more that $\partial_{x}^{1/2}$ smoothness while the one for the KP-II
gains the full derivative $\partial_{x}$.

The first result regarding well-posedness for a KP type equation with
low regularity is due to Ukai \cite{U}. He uses a standard energy
method that does not recognize the type I or II of the equation. His
result provides local well-posedness for initial data and their
antiderivatives  in $\in H^{s}, \, s\geq 3$.
Faminskii \cite{F}
observed a better smoothing effect in the KP
II evolution and used this to prove well-posedness results.
Bourgain performed a Fourier analysis \cite{B} 
of the term $\partial_x u^2$
in the KP-II equation in which the derivative is recovered in a nonlinear
way. The result obtained gave local well-posedness
of KP-II for initial data in $L^2 $. Since the $L^2$ norm is conserved
during the KP-II evolution, the $L^2$ local result may be iterated to
prove global well-posedness.
Takaoka \cite{Ta1} and Takaoka and Tzvetkov \cite{TT}
improved Bourgain's result by proving local well-posedness in an
anisotropic Sobolev space  $H^{-1/3+\epsilon,0}_{x,y}$.
For the KP-I equation the situation is more delicate. There are
several results on local and global existence of solutions, but not a
satisfactory well-posedness theory for data with no more than
two derivatives. Fokas and Sung \cite{FS}, and Zhou \cite{Z}, obtained  global
existence  for small data via inverse scattering
techniques.
Schwarz \cite{Sw} proved existence of weak global
periodic solutions with small $L^{2}$ data. The smallness condition was
subsequently removed \cite{CKPI}. Tom \cite{T} 
proved existence of global weak solutions
for initial data in $H^{1}$ together with their antiderivative.
For well-posedness results, we recall the work
of Saut \cite{S}, Isaza, Mej\'ia and Stallbohom~\cite{IMS} and
finally
the work of I\'orio and Nunes~\cite{IN}. The last two authors
use the quasi linear theory of Kato, together with parabolic
regularization, to prove local well-posedness with data and their
antiderivatives in $H^{s}, s>2$. The limitation $s>2$  is
needed in order to insure that $\partial_{x}u\in L^{\infty}$, an
essential assumption for the proof. Molinet, Saut and Tzetkov
\cite{MST1}
also proved that if one is willing to assume
more regularity for the initial data (at least three derivatives in the
$x$ variable and two in the $y$ variable need to be in $L^{2}$),
then global well-posedness holds.

In this paper we use a method involving oscillatory integrals to
prove that for small
\footnote{The precise ``smallness'' condition can be found in
Theorem \ref{main1} and \ref{main2}}
initial data $u_{0}$ in a certain weighted Sobolev
space
\footnote{The precise definition of this space can be found in
(\ref{id}).}, defined with at most two derivatives,
the IVP (\ref{ivp1}) is globally well-posed. The fact that we had to
use ``weights'' in the definition of our space agrees with some recent
counterexamples of Molinet, Saut and Tzetkov \cite{MST}. These
counterexamples suggest
``that any iterative method applied to the integral formulation of
the KP-I equation always fails'' when the initial data are
only in anisotropic Sobolev spaces.

The method of proof that we adopt here follows the approach used by
Kenig, Ponce and Vega to treat the Schr\"odinger IVP with derivative
in the non-linearity \cite{KPV3}. Here the situation is more complex
due to the anisotropic  nature of the problem. A weaker version 
of the so called ``smoothing effect estimates'' and
``maximal function estimates'', appeared in the work of
Isaza-Mej\'ia-Stallbohm~\cite{IMS1},
but these estimates were not strong enough to complete  a
fixed point argument.

In the rest of this section, we introduce some notation and
definitions. Then, in Section 2, we state the main theorem. In
Section 3, we present the estimates related to the maximal function
associated to (\ref{lsolution}).  Section 4 is dedicated to the smoothing
effect estimates and Section 5 to the group estimates. Finally, in
Section 6, we present the main steps of the proof of the
well-posedness theorem
via the fixed point argument. The paper has also a short appendix on
some consequences of the fractional Leibniz rule.

\vspace{1cm}
\noindent
{\bf{ Notation:}} We denote the Fourier transform of a function
$f(x,y)$ by
\[\F(f)(\x,\lambda)=\hat{f}(\x,\lambda)=\int_{\dbR^{2}}
e^{i(\x x+\lambda y)}f(x,y)dx dy\]
and the inverse Fourier transform by
\[\F^{-1}(g)(\x,\lambda)=\check{g}(x,y)=\int_{\dbR^{2}}
e^{-i(\x x+\lambda y)}g(\x,\lambda)d\x d\lambda.\]
(We systematically ignore various ``2 $\pi$-constants''.)
We will often use the symbol $\langle y\rangle=(1+y^{2})^{1/2}$
and the operators $D_{x}^{\sigma}$ and $D_{y}^{\gamma}$ which
are defined through the Fourier transform as the multiplier operators
$\F(D_{y}^{\gamma}f)(\x,\lambda)=
|\lambda|^{\gamma}\hat{f}(\x,\lambda)$ and
$\F(D_{x}^{\sigma}f)(\x,\lambda)=
|\x|^{\sigma}\hat{f}(\x,\lambda)$.
We denote with
$H^{\dot{\sigma},\dot{\gamma}}_{\langle y\rangle^\alpha}$, the
closure of the Schwartz functions with respect to the norm
\begin{equation}
\|f\|_{H^{\dot{\sigma},\dot{\gamma}}_{<y>^\alpha}}=
\|<y>^\alpha D_x^\sigma D_y^{\gamma} f\|_{L^2_{\dbR^2}}.
\label{hspace}\end{equation}
We remove the dot on the indices $\sigma$ and $\gamma$
if we substitute $D_x^\sigma$ with $(1+D_x^\sigma)$ and
$D_y^\gamma$ with $(1+D_y^\gamma)$ respectively.
We will also use a variety of mixed $L^{p}$ norms. For example
the space $L^{r}_{t}L^{q}_{x}L^{p}_{y}$ is the space of functions
equipped with the norm
$$
\|f\|_{L^{r}_{t}L^{q}_{x}L^{p}_{y}}=
\left(\int\left(\int
\left(\int |f|^{p}(x,y,t) \,dy\right)^{q/p} dx\right)^{r/q}dt
\right)^{1/r}.$$
We will also write $L^{r}_{t}L^{p}_{x}L^{p}_{y}=
L^{r}_{t}L^{p}_{xy}$ and $L^{2}_{g(x,y)dxdy}$ to indicate
the space of $L^{2}$ functions with respect to the measure
$g(x,y)dxdy$.

Next, we define some projection
operators  that will appear throughout the paper.
The operator
$P_{+}:L^2(\dbR^2)\longrightarrow L^2(\dbR^2)$ is defined as
\begin{equation}
\widehat{P_{+}u}(\xi,\lambda)=\chi_{\{|\xi|\geq
\left|\frac{\lambda}{\x}\right|\}}\hat{u}(\xi,\lambda).
\label{P+}\end{equation}
and $P_{-}=Id-P_{+}$.
It will become clear later that the estimates on the solution of the
linear problem associated to (\ref{ivp}) will be easier to obtain if
one could assume that all  frequencies $\x$ were far from zero. To
put ourselves in this setting,  we will use the projection
operator  $Q:L^2(\dbR^2)\longrightarrow L^2(\dbR^2)$ defined by
the formula
\begin{equation}
\widehat{Qu}(\xi,\lambda)=\chi_{\{|\xi|\geq 1\}}\hat{u}(\xi,\lambda).
\label{q}\end{equation}
We use the notation $A\lesssim B$ to
indicate that there exists a constant $m\ne 0$ such that
$A\leq mB$. We will use the abbreviations L.H.S. and R.H.S. to refer
to the left hand side and right hand side of inequalities or
equations.


\section{The Main Theorem}

We now consider the KP-I initial value problem

\begin{equation}
\left\{ \begin{array}{l}
\partial_t u + \partial_x^3 u - \partial_x^{-1}\partial_y^2u
+\beta\partial_xu^2= 0, \\
u(x,0) = u_0(x), \hspace{1.5cm}(x,y) \in \dbR^2, \,  t \in \dbR,
\end{array}\right.
\label{ivp}\end{equation}
where $u=u(t,x,y)$ is a scalar unknown function and $\beta$
is a real constant.

Using the operator $Q$ defined at the end of Section 1, which
basically selects the large frequencies,
we define the space\footnote{See also Definition \ref{spacexy0}.}
$Z_{0}$. We say that  $u_{0} \in Z_{0}$ if
\begin{eqnarray}
  \label{id} \|u_{0}\|_{Z_{0}}&=&\|u_{0}\|_{H^{2}}\\
\nonumber  &+&
   \|(Id-Q)u_{0}\|_{H^{-\dot{\sigma_{0}},\gamma_{0}}_{\langle
   y\rangle^{\alpha}}}+
   \|(Id-Q)u_{0}\|_{H^{-\dot{\sigma_{1}},\gamma_{1}}_{\langle
   y\rangle^{\alpha}}}+
   \|(Id-Q)u_{0}\|_{H^{-\dot{\sigma_{2}},\gamma_{2}}}\\
   \nonumber&+&
   \|Qu_{0}\|_{H^{\sigma_{3},\gamma_{3}}_{\langle
   y\rangle^{\alpha}}}+
   \|Qu_{0}\|_{H^{0,\gamma_{4}}_{\langle
   y\rangle^{\alpha}}}+
   \|Qu_{0}\|_{H^{\sigma_{4},0}_{\langle
   y\rangle^{\alpha}}}<\infty,
\end{eqnarray}
where $\alpha=\gamma_{0}=\sigma_{0}=\gamma_{3}=1/2+\epsilon,
\sigma_{1}=1/4+\epsilon, \sigma_{2}=\sigma_{3}=3/4+\epsilon,
\sigma_{4}=\gamma_{1}=\gamma_{4}=1+\epsilon,
\gamma_{2}=3/2+\epsilon,$ where $\epsilon>0$
and small.

The first result we present is a
global well-posedness statement for small initial data:
\begin{theorem}
For any $T>0$, there exists $\delta>0$ such that for any
$u_0 \in Z_{0}$, and
\begin{equation}
    \max(\|\langle y\rangle^{\alpha}
Q u_0\|_{L^{2}}, \|\langle y\rangle^{\alpha}D_{x}^{-\sigma_{0}}
(Id-Q) u_0\|_{L^{2}},
\|D_{x}^{-\sigma_{3}}(Id - Q) Qu_0\|_{L^{2}})\leq \delta,
\label{smallness}\end{equation}
there exists a unique solution
$u(x,y,t)$ for the IVP (\ref{ivp})
in the interval $[0,T]$, satisfying
\[u\in C([0,T],Z_{0})\cap Z_{T},\]
where the space $Z_{T}$ is defined in Definition \ref{spacexy}.
Moreover, for any $T'\in (0,T)$, there exists $\rho>0$ such that
the map $\widetilde{u_0}\longrightarrow \tilde{u}$ from
$\{\widetilde{u_0}\in Z_{0}/\|\widetilde{u_0}-u_0\|_
{Z_{0}}\leq \rho\}$ into $C([0,T'], Z_{0})\cap Z_{T'}$ is Lipshitz.
\label{main1}\end{theorem}
In this second theorem, we relax \eqref{smallness} a little, but we
only obtain a local well-posedness result:
\begin{theorem}
There exists $\delta>0$ such that for any $u_0 \in Z_{0}$, and
$$\|\langle y\rangle^{\alpha}(1+D_{y})^{\gamma_{4}}
u_0\|_{L^{2}}\leq \delta,$$
there exists $T=T(\|u_{0}\|_{Z_{0}})$
and a unique solution $u(x,y,t)$ for the IVP (\ref{ivp})
in the interval $[0,T]$, satisfying
\[u\in C([0,T],Z_{0})\cap Z_{T}.\]
Moreover, for any $T'\in (0,T)$, there exists $\rho>0$ such that
the map $\widetilde{u_0}\longrightarrow \tilde{u}$ from
$\{\widetilde{u_0}\in Z_{0}/\|\widetilde{u_0}-u_0\|_
{Z_{0}}\leq \rho\}$ into $C([0,T'], Z_{0})\cap Z_{T'}$ is Lipshitz.
\label{main2}\end{theorem}
Using the operator $Q$ again,
we transform (\ref{ivp}) into
a system with unknowns  $Qu=u_1$ and $(Id-Q)u=u_2$:

\begin{equation}
\left\{ \begin{array}{l}
\partial_t u_1 + \partial_x^3 u_1 - \partial_x^{-1}
\partial_y^2u_1
+\beta Q\partial_x(u_1^2+u_{2}^{2}+2u_{1}u_{2})= 0 \\
\partial_t u_2 + \partial_x^3 u_2 - \partial_x^{-1}
\partial_y^2u_2
+\beta(Id-Q)\partial_x(u_1^2+u_2^2+2u_2u_1)= 0 \\
u_1(x,0) = Qu_0(x)=w_{0}\\
u_2(x,0) = (Id-Q)u_0(x)=v_{0},  \hspace{1.5cm}(x,y)
\in \dbR^2, \,  t \in \dbR.
\end{array}\right.
\label{ivp2}\end{equation}
Now observe that for any function $f, \,(Id-Q)\partial_x f$
and $ f$, have similarly sized Fourier transforms supported on 
low frequencies. Hence, the effective nonlinear term for the problem 
is concentrated in the first equation,
that is the term $Q\partial_x(u_{1}^{2}+u_2^2+2u_2u_1)$. On the
other hand the dispersive function for the first equation is
$\phi(\xi,\lambda)=\chi_{\{|\xi|>1\}}(\xi^3+\lambda^2/\xi)$,
which is not singular. 

\section{The Maximal Function Estimate}
In this section we prove a maximal function type estimate for the
solution of the linear initial value problem associated to (\ref{ivp}).
Consider the problem
\begin{equation}
\left\{ \begin{array}{l}
\partial_t u + \partial_x^3 u - \partial_x^{-1}
\partial_y^2u= 0, \\
u(x,0) = u_0(x) \hspace{1.5cm}(x,y) \in \dbR^2, \,  t \in \dbR,
\end{array}\right.
\label{livp}\end{equation}
and denote with $U(t)u_{0}$  its solution. It is easy to see that
\begin{equation}
 U(t)u_{0}(x)=\int_{\dbR^{2}}e^{i(t(\x^{3}+\frac{\lambda^{2}}{\x})
 +\x x+\lambda y)} \widehat{u_{0}}(\x,\lambda)d\x d\lambda.
\label{lsol}\end{equation}
We have the following theorems:
\begin{theorem}
For any $\sigma>3/4, \gamma>1/2$ and $\theta>1$,
\begin{eqnarray}
\label{tmaxx}\left(\sum_{s=-\infty}^{\infty}\sup_{ s\leq y<s+1}
\sup_{|t|\leq 1}\sup_{x}|QU(t)u_0(x,y)|^2\right)^{1/2}
&\lesssim& \|(1+D_x)^{\sigma}(1+D_y)^{\gamma}Qu_0\|_{L^2(\dbR^2)}\\
\label{tmaxy}\left(\sum_{r=-\infty}^{\infty}\sup_{ r\leq x<r+1}
\sup_{|t|\leq 1}\sup_{y}|QU(t)u_0(x,y)|^2\right)^{1/2}
&\lesssim& (\|(1+D_x)^{\sigma}(1+D_y)^{\gamma}Qu_0\|_{L^2(\dbR^2)}\\
\nonumber &+&\|(1+D_y)^{\theta}Qu_0\|_{L^2(\dbR^2)}).
\end{eqnarray}
\label{max}\end{theorem}
For the low frequency solution we have:
\begin{theorem}
For $\sigma>1/4, \gamma>1/2$ and $\theta>1$,
\begin{eqnarray}
\label{tmaxxlow}& &\left(\sum_{s=-\infty}^{\infty}\sup_{ s\leq y<+1}
\sup_{|t|\leq 1}\sup_{x}|(Id-Q)U(t)u_0(x,y)|^2\right)^{1/2}\\
\nonumber&\lesssim& \|D_x^{-\gamma}(1+D_y)^{\gamma}(Id-Q)u_0\|_{L^2(\dbR^2)}\\
\label{tmaxylow}& &\left(\sum_{r=-\infty}^{\infty}\sup_{ r\leq x<r+1}
\sup_{|t|\leq 1}\sup_{y}|(Id-Q)U(t)u_0(x,y)|^2\right)^{1/2}\\
\nonumber&\lesssim& 
\|D_{x}^{-\sigma}(1+D_y)^{\theta}(Id-Q)u_0\|_{L^2(\dbR^2)}).
\end{eqnarray}
\label{maxlow}\end{theorem}
The proof of these  theorems follows arguments presented 
in \cite{KPV3}. We use the following lemma:
\begin{lemma}
Let $\phi(\xi,\lambda)=\xi^3+\lambda^2/\xi$ and let $\theta_{k,j}$
be a $C^{\infty}$ function supported on the set $\Delta_{k,j}=
\{|\xi|\sim 2^k \mbox{ and } |\lambda/\xi|\sim 2^j, \, \, \mbox{ for }
k\in \dbZ, j\in \dbZ\}$. Then the function
\begin{equation}
I_{k,j}^t(x,y)=\int_{\dbR^2}e^{i(t\phi(\xi,\lambda)+(x,y)(\xi,\lambda)}
\theta_{k,j}(\xi,\lambda)\, d\xi d\lambda
\end{equation}
satisfies

\[|I_{k,j}^t(x,y)|\leq H_{k,j}(|x|,|y|),\, \, \,
\mbox{ for $|t|\leq 1$}\]
and
\begin{eqnarray}
\label{maxx}\sum_{-\infty<s<\infty}\sup_{r}H_{k,j}(|r|,|s|)
&\lesssim&  2^{\alpha(k,j)},\\
\label{maxy}\sum_{-\infty<r<\infty}\sup_{s}H_{k,j}(|r|,|s|)
&\lesssim&  2^{\delta(k,j)},
\end{eqnarray}
where $\alpha(k,j)=\frac{5}{2}k+j$ if $k\geq 0$ and
$\alpha(k,j)= j+$ if $k\leq 0$, and
\begin{equation*}
  \delta(k, j) = 
\left\{
\begin{matrix}
\frac{5k}{2} + j  &~{\mbox{if}}~ k \geq \max (0,j) ,\\
\frac{3k}{2} + 2j & ~{\mbox{if}}~ j \geq k \geq 0 ~{\mbox{or}}~ [ k < 0 ~{\mbox{and}}~ k \geq 2j], \\
2k + j & if k  < \min ( 0, 2j ).
\end{matrix} 
\right.
\end{equation*}

\label{I}\end{lemma}

\begin{proof}
We subdivide the interval of time $[0,1]$ into dyadic intervals of
size $[2^{\sigma-1}, 2^{\sigma}]$, for $\sigma\leq 0$. Then we 
rescale the integral
$I_{k,j}$ by setting $t=\rho 2^{\sigma}, \, \rho\in [1/2,1]$ and
\begin{equation}
2^{\sigma/3}\x=\z \, \, \, \mbox{ and } \, \, \,
2^{2\sigma/3}\la=\mu.
   \end{equation}
Then
\begin{equation}
I_{k,j}(x,y)= \int_{\dbR^2}
2^{-\sigma}e^{i(\rho\phi(\z,\mu)+(2^{-\sigma/3}x,
2^{-2\sigma/3}y)(\z,\mu)}
\theta_{\tilde{k},\tilde{j}}(\z,\mu)\, d\z d\mu =
2^{-\sigma}I_{\tilde{k},\tilde{j}}(\tilde{x},\tilde{y})
\label{in}\end{equation}
where
\begin{eqnarray}
\label{tildekj}    \tilde{k} = k+\sigma/3, & &
\tilde{j} = j+\sigma/3,\\
\label{tildexy}\tilde{x}=2^{-\sigma/3}x, & &
\tilde{y}=2^{-2\sigma/3}y.
\end{eqnarray}
For simplicity we drop the ``tilde'' from $k, j, x $  and $y$
and we keep in mind that again $k, j\in \dbZ$. We start by proving
(\ref{maxy}). We make the decomposition
\begin{eqnarray*}
   I_{k,j}(x,y)&=&I_{k,j}(x,y)\chi_{\{|x|\leq 1\}}+
   I_{k,j}(x,y)\chi_{\{|x|\geq C\max{(2^{2j},2^{2k})}\}}+
I_{k,j}(x,y)\chi_{\{|x|\leq C\max{(2^{2j},2^{2k})}\}}
\\
&=&K_{k,j}^{0}(x,y)+K_{k,j}^{1}(x,y)+K_{k,j}^{2}(x,y).
\end{eqnarray*}
for an appropriate $C \gg 1$ to be determined later. A different 
decomposition will be used later when we address \eqref{maxy}.

\vskip .2in
\noindent
{\bf{ Estimate of $K^0_{k,j}$}}\\
Here we simply have 
\begin{equation}
|K_{k,j}^{0}(x,y)|\leq\chi_{\{|x|\leq 1\}}\int_{\dbR^2}
\theta_{k,j}(\z,\mu)\,d\z d\mu=2^{2k+j}\chi_{\{|x|\leq 1\}}
=H^{0}_{k,j}(|x|,|y|).
\label{k0}\end{equation}
Define the function
\[\Psi_{x,y}(\z,\mu)=\rho\phi(\z,\mu)+(x,y)\cdot(\z,\mu)\]
for $\rho \in [1/2,1]$. In the rest of the proof we will drop $\rho$.
Compute
\begin{equation}
   \nabla\Psi_{x,y}(\z,\mu)=[3\z^2-\mu^2/\xi^2+x, 2\mu/\z+y]
\label{grad}\end{equation}
and
\begin{eqnarray}
\label{partialmu}\partial_\mu^2\Psi= 2/\z & &
\partial_\mu\theta_{k,j}\sim
2^{-k-j}\tilde\theta_{k,j}\\
\label{partialz}\partial_\z^2\Psi= 6\z+\mu^2/\z^3 & &
\partial_\z\theta_{k,j}\sim
2^{-k}\tilde\theta_{k,j}.
\end{eqnarray}
We now choose the constant $C$ in the defining condition of $K^i_{k,j},~j =
0,1,2,$ so that, in the present region under consideration $K^1_{k,j}$, we have
\[|\partial_{\z}\Psi_{x,y}(\z,\mu)|\gtrsim |x|.\]
%

\vskip .2in
\noindent
{\bf{ Estimate of $K^1_{k,j}$}}\\
Observe that in this case we can assume that  $\max(k,j)\geq 0$. (Otherwise
we are in the region where $|x| \leq 1$ which was considered above.)
We integrate by parts with respect to $\z$ twice and we use
(\ref{partialz}) to obtain
\begin{eqnarray}
\nonumber |K_{k,j}^{1}(x,y)|&\lesssim&\int
\left(\frac{|\theta_{\z\z}|}{|\Psi_{\z}|^{2}}+
\frac{|\theta_{\z}\Psi_{\z\z}|}{|\Psi_{\z}|^{3}}+
\frac{|\theta\Psi_{\z\z\z}|}{|\Psi_{\z}|^{3}}+
\frac{|\theta\Psi_{\z\z}^{2}|}{|\Psi_{\z}|^{4}}\right)d\z d\mu\\
\nonumber &\lesssim& 2^{2k+j}\left(\frac{2^{2k}}{\max{(|x|^{2},1)}}+
\frac{2^{k}\max{(2^{k},2^{2j-k})}}{\max{(|x|^{3},1)}}+
\frac{\max{(1,2^{2j-2k})}}{\max{(|x|^{3},1)}}+
\frac{\max{(2^{2k},2^{4j-2k})}}{\max{(|x|^{4},1)}}\right)\\
\label{k1}&\lesssim& \frac{2^{2k+j+}}{|x|^{1+}}
\chi_{\{|x|>C\max{(1, 2^{2j},2^{2k})}\}}=H^{1}_{k,j}(|x|,|y|),
\end{eqnarray}
%

\vskip .2in
\noindent
{\bf Estimate of $K^{2}_{k,j}$}\\
In this case $\partial_{\z}\Psi$ could be zero and no integration by
parts could be performed. We need the following lemma
\begin{lemma}
Assume $\mu$ is fixed and that there exists $\z_0= \z_{0}(\mu)
\in \Delta_{k,j}$ such that $\partial_{\z}\Psi(\z_0,\mu)=0$. Then
\begin{equation}
|I_{k,j}(x,y)|=|\int e^{i(\z^3+\mu/\z-(x,y)(\z,\mu))}
\theta_{k,j}
(\z,\mu)\,d\mu d\z|\lesssim 2^{\beta(k,j)},
\end{equation}
where $\beta(k,j)=\frac{k}{2}+j$ if $k\geq j$ and
$\beta(k,j)=\frac{3k}{2}$ if $k< j$.
\label{vdc}\end{lemma}
\begin{proof}
We use the Van der Corput lemma (see for example Corollary of
Proposition 2 Chap. VIII, in \cite{Stein}). We first make
a change of variables so that $\Delta_{k,j}$ is
transformed into $\Delta_{1,1}$, that is we set $(\xi,\tau)=
(2^{-k}\z, 2^{-j}\mu)$. With the new variables we have
$$\Psi(\xi,\tau)=2^{3k}\x^3+2^{2j-k}\tau^2/\x+(2^{k}x,2^{j}y)
\cdot(\x,\tau)=m_{k,j}\tilde{\Psi}(\xi,\tau)$$
where $m_{k,j}=\max{(2^{3k},2^{2j-k})}$. Then
\begin{equation}
   |K_{k,j}^{2}(x,y)|\sim 2^{2k+j}
   \left|\int e^{im_{k,j}\tilde{\Psi}(\xi,\tau)}\theta_{1,1}
(\xi,\tau)\,d\tau d\xi\right|.
\end{equation}
It is  easy to check that $\|\tilde{\Psi}\|_{C^{3}}\leq C$ and that
\begin{equation}
|\partial_{\x\x}\tilde{\Psi}|=m_{k,j}^{-1}|\xi|\left(2^{3k}6+2
2^{2j-k}\frac{\tau^{2}}{\xi^{4}}\right)\sim 1,
\end{equation}
hence by integrating first with respect to $\x$ using Van der Corput
lemma, and then with respect to $\tau$, we obtain
\begin{equation}
|K_{k,j}^{2}(x,y)|\sim 2^{2k+j}(m_{k,j})^{-1/2}=2^{\beta(k,j)}
\chi_{\{1\leq |x|\lesssim \max{(1,2^{2j},2^{2k})}\}}=H^{2}(|x|,|y|)
\label{k2}\end{equation}
where $\beta(k,j)=k/2+j$ if $k\geq j$ and $\beta(k,j)=3k/2$ if $k<j$
and $\max{(k,j)}>0$.
\end{proof}
We now go back to the ``tilde'' notation and we
define
\[\tilde{H}_{\tilde{k},\tilde{j}}(\tilde{x},\tilde{y})
=\sum_{i=0,1,2}\tilde{H}_{\tilde{k},\tilde{j}}^{i}
(\tilde{x},\tilde{y}).\]
Using  (\ref{tildexy}), (\ref{k0}), (\ref{k1}) and (\ref{k2}) we
finally have that
\begin{eqnarray*}
   \sum_{r\in\dbZ }\sup_{s}H_{k,j}(|r|,|s|)&=&
2^{-\sigma}\sum_{r\in\dbZ }\sup_{s}\tilde{H}_{\tilde{k},\tilde{j}}
(2^{-\sigma/3}|r|,2^{-2\sigma/3}|s|)\\
&=&2^{-2\sigma/3}\sum_{\bar{r}\in\dbZ }
\sup_{\bar{s}}H_{\tilde{k},\tilde{j}}(|\bar{r}|,|\bar{s}|)\\
&\leq&2^{-2\sigma/3}2^{\tilde\delta(\tilde{k},\tilde{j})}
\end{eqnarray*}
where
\begin{equation}
\tilde\delta(\tilde{k},\tilde{j})=\left\{\begin{array}{l}
\frac{5\tilde{k}}{2} +\tilde{j}  \, \, \, \mbox{ if }\, \,
\tilde{k}\geq \max{(0,\tilde{j})}\\
\frac{3\tilde{k}}{2} +2\tilde{j}  \, \, \, \mbox{ if }\, \,
\tilde{j}\geq \max{(0,\tilde{k})}\\
2\tilde{k}+\tilde{j}\, \, \, \mbox{ if }\, \, \tilde{k},\tilde{j}
\leq 0,
\end{array}\right.
\end{equation}
which by (\ref{tildekj}) gives
\[\sum_{r\in\dbZ }\sup_{s}H_{k,j}(|r|,|s|)\lesssim
2^{\sigma/3+\delta(k,j)}, \, \, \, \sigma\leq 0\]
where
\begin{equation}
\delta(k,j)=\left\{\begin{array}{l}
\frac{5k}{2} +j  \, \, \, \mbox{ if }\, \,
k\geq \max{(0,j)}\\
\frac{3k}{2} +2j  \, \, \, \mbox{ if }\, \,
j\geq k\geq 0 \, \, \mbox{ or }
k<0 \, \, \mbox{ and } \, \, k\geq 2j\\
2k+j\, \, \, \mbox{ if }\, \, k\leq \min{(0,2j)}.
\end{array}\right.
\end{equation}
This proves (\ref{maxy}).  We now use similar ideas to prove
(\ref{maxx}). After we rescaled the time as we did above we write

\begin{eqnarray*}
   I_{k,j}(x,y)&=&I_{k,j}(x,y)\chi_{\{|y|\leq 1\}}+
I_{k,j}(x,y)\chi_{\{|y|\geq \max{(1,C2^j)}\}}
+I_{k,j}(x,y)\chi_{\{1\leq |y|\lesssim 2^j\}}\\
&=&K_{k,j}^{3}(x,y)+K_{k,j}^{4}(x,y)+K_{k,j}^{5}(x,y).
\end{eqnarray*}

\vskip .2in
\noindent
{\bf{ Estimate of $K^{3}_{k,j}$}}\\
Clearly we can estimate
\begin{equation}
|K_{k,j}^{3}|(x,y)\lesssim 2^{2k+j}\chi_{\{|y|\leq 1\}}=H^{3}_{k,j}.
\label{k3}\end{equation}

\vskip .2in
\noindent
{\bf{ Estimate of $K^4_{k,j}$}}\\
Observe first that in this case we can assume that  $j\geq 0$.
We then use the fact that for $C\gg 1$
\[|\partial_{\mu}\Psi|\gtrsim |y|,\]
and we can integrate by parts with respect to $\mu$ twice and  use
(\ref{partialmu}) to obtain
\begin{eqnarray}
\nonumber |K_{k,j}^{4}(x,y)|&\lesssim&\int
\left(\frac{|\theta_{\mu\mu}|}{|\Psi_{\mu}|^{2}}+
\frac{|\theta_{\mu}\Psi_{\mu\mu}|}{|\Psi_{\mu}|^{3}}+
\frac{|\theta\Psi_{\mu\mu\mu}|}{|\Psi_{\mu}|^{3}}+
\frac{|\theta\Psi_{\mu\mu}^{2}|}{|\Psi_{\mu}|^{4}}\right)d\z d\mu\\
\label{k4}&\lesssim& \frac{2^{-2j+}}{|y|^{1+}}
\chi_{\{|y|>C\max{(1, 2^{j})}\}}=H^{4}_{k,j}(|x|,|y|),
\end{eqnarray}

\vskip .2in
\noindent
{\bf{ Estimate of $K^{5}_{k,j}$}}\\
In this case $\partial_{\mu}\Psi$ could be zero and no integration by
parts could be performed. Instead we use an argument similar to the one
presented during the proof of Lemma \ref{vdc}. Rescale the function
$\Psi$ using variables $\z$ and $\tau$ and define $m_{k,j}=2^{2j-k}$.
Then
\begin{equation}
|K_{k,j}^{5}(x,y)|\sim 2^{2k+j}(m_{k,j})^{-1/2}=2^{3k/2}
\chi_{\{1<|x|\lesssim 2^{j}\}}=H^{5}(|x|,|y|).
\label{k5}\end{equation}
using  (\ref{tildexy}), (\ref{k3}), (\ref{k4}) and (\ref{k5}) we
finally have that
\begin{eqnarray*}
   \sum_{s\in\dbZ }\sup_{r}H_{k,j}(|r|,|s|)&=&
2^{-\sigma}\sum_{r\in\dbZ }\sup_{s}\tilde{H}_{\tilde{k},\tilde{j}}
(2^{-\sigma/3}|r|,2^{-2s/3}|s|)\\
&=&2^{-\sigma/3}\sum_{\bar{r}\in\dbZ }
\sup_{\bar{s}}H_{\tilde{k},\tilde{j}}(|\bar{r}|,|\bar{s}|)\\
&\leq&2^{-\sigma/3}2^{\tilde\delta(\tilde{k},\tilde{j})}
\end{eqnarray*}
where
\begin{equation}
\label{missingref}
\tilde\delta(\tilde{k},\tilde{j})=
\max{(2\tilde{k}+\tilde{j},-\chi_{\{\tilde{j}\geq 0\}}2\tilde{j},
\chi_{\{\tilde{j}\geq 0\}}(5/2\tilde{k}+\tilde{j}))}.
\end{equation}
After we
replace $\tilde k$ with $k$ and $\tilde j$ with $j$, we have
\[\sum_{s\in\dbZ }\sup_{r}H_{k,j}(|r|,|s|)\leq 2^{-\sigma\epsilon}
2^{\alpha(k,j)}\]
where $\alpha(k,j)=\frac{5}{2}k+j+$ if $k\geq 0$ and
$\alpha(k,j)=j+$ if $k\leq  0$
and this proves (\ref{maxx}).
\end{proof}

\noindent
The proofs of Theorems \ref{max} and \ref{maxlow} follow the
basic steps of the proof of Theorem 3.2 in \cite{KPV3}.
\begin{proof}[Proof of Theorems \ref{max} and \ref{maxlow}:] 
Here we prove only  (\ref{tmaxx}), (\ref{tmaxy}) will follow from
similar arguments and Lemma \ref{I}. Let $\theta_{k,j}$ be as in
in the proof of Lemma \ref{I}. We define
\begin{equation}
\F(U_{k,j}(t)u_0)(\xi,\lambda)=e^{it\phi(\xi,\lambda)}
\theta_{k,j}(\xi,\lambda)\widehat{u_0}(\xi,\lambda)
\end{equation}
with $k\in \dbN$ and $j\in \dbZ$. It suffices to show that
\begin{equation}
\left(\sum_{s=-\infty}^{\infty}\sup_{s\leq y<s+1}
\sup_{|t|\leq 1}\sup_{x}
|U_{k,j}(t)u_0(x,y)|^2\right)^{1/2}
\lesssim 2^{\alpha(k,j)}\|u_0\|_{L^2(\dbR^2)}
\label{dyadicmaxx}\end{equation}
where $\alpha(k,j)$ is defined in Lemma \ref{I}.
We recall that the dual operator to
$U_{k,j}(t)$, with $t\in [-1,1]$,
is the operator $T^{*}$ such that
$$T^{*}g(t,x,y)=\int_{-1}^{1}U_{k,j}(t)g(t,x,y) dt.$$
By duality is it enough to show that
\begin{equation}
\left\|\int_{-1}^{1}U_{k,j}(t)g(t,\cdot)\,dt\right
\|_{L^2}\lesssim 2^{\alpha(k,j)}
\left(\sum_{s=-\infty}^{\infty}\left(\int_s^{s+1}\int_{-1}^1\int_{\dbR}
|g(x,y,t)|\,dx\,dt\,dy\right)^{2}\right)^{1/2}.
\label{dual}\end{equation}
Using P. Tomas'argument on L.H.S of (\ref{dual}), it follows that
\begin{eqnarray*}
& &\left\|\int_{-1}^{1}U_{k,j}(t)g(t,\cdot)\,dt\right
\|_{L^2}^{2}=\int_{\dbR^{2}}\int_{-1}^{1}
\left(\int_{-1}^{1}U_{k,j}(t-\tau)g(\tau,x,y) d\tau\right)
g(t,x,y) \,dt \,dx \,dy\\
&\lesssim& \sum_{s}\sup_{s\leq
y<s+1}\int_{s}^{s+1}\int_{\dbR}\int_{-1}^{1}
\left(\int_{-1}^{1}U_{k,j}(t-\tau)g(\tau,x,y) d\tau\right)
g(t,x,y) \,dt \,dx \,dy\\
&\lesssim& \sum_{s}\left(\sup_{s\leq y<s+1}\sup_{x}\sup_{|t|\leq 1}
\left|\int_{-1}^{1}U_{k,j}(t-\tau)g(\tau,x,y)\,d\tau\right|\right)
\left(\int_{s}^{s+1}\int_{\dbR}\int_{-1}^{1}|g(t,x,y)|\,dt \,dx \,dy
\right)\\
&\lesssim& \left(\sum_{s}\left(\sup_{s\leq
y<s+1}\sup_{x}\sup_{|t|\leq 1}
\left|\int_{-1}^{1}U_{k,j}(t-\tau)g(\tau,x,y)\,d\tau
\right|\right)^{2}\right)^{1/2}\\
&\times&\left(\sum_{s}\left(\int_{s}^{s+1}\int_{\dbR}\int_{-1}^{1}|
g(t,x,y)|\,dt \,dx \,dy\right)^{2}\right)^{1/2}.
\end{eqnarray*}
So in order to prove (\ref{dyadicmaxx}) it is enough to show that
\begin{eqnarray}
\nonumber& &\left(\sum_{s=-\infty}^{\infty}
\sup_{s\leq y<s+1}\sup_{x}\sup_{|t|\leq 1}
\left|\int_{-1}^1U_{k,j}(t-\tau)g(\tau,x,y)\,d\tau\right|^2
\right)^{1/2}\\
&\lesssim&2^{2\alpha(k,j)}\left(\sum_{s=-\infty}^{\infty}
\left(\int_s^{s+1}\int_{-1}^1\int_{\dbR}
|g(x,y,t)|\,dx\,dy\,dt\right)^2\right)^{1/2}.
\end{eqnarray}
We observe that
$$
\left|\int_{-1}^1U_{k,j}(t-\tau)g(\tau,x,y)\,d\tau\right|\lesssim
\int H_{k,j}(z,w)\int_{-1}^1|g(x-z,y-w,\tau)|\,dz\,dw\,d\tau
$$
and by Young's inequality,
$$
\sup_{x}\left|\int_{-1}^1U_{k,j}(t-\tau)g(\tau,x,y)\,d\tau\right|
\lesssim
\int \sup_{z}H_{k,j}(z,w)
\left(\int_{\dbR}\int_{-1}^1|g(z,y-w,\tau)|
\,d\tau\,dz\right)\,dw.
$$
If we partition $w$ we can continue with
$$
\lesssim\sum_{\bar s=-\infty}^\infty \sup_{z}
\sup_{\bar{s}\leq w<\bar{s}+1}H_{k,j}(|z|,|w|)
\int_{\dbR}\int_{\bar{s}}^{\bar{s}+1}\int_{-1}^1|g(z,y-w,\tau)|
d\tau\,dw\, dz.
$$
Then we partition $y$ with $s$ and take the $L^{2}$ norm,
to find
\begin{eqnarray*}
& \lesssim & \left(\sum_{s=-\infty}^{\infty}\left(\sup_{s\leq y<s+1}
\sum_{\bar s=-\infty}^\infty \sup_{z}
\sup_{\bar{s}\leq w<\bar{s}+1}H_{k,j}(|z|,|w|)
\int_{\dbR}\int_{\bar{s}}^{\bar{s}+1}\int_{-1}^1|g(z,y-w,\tau)|
d\tau\,dw\, dz\right)^{2}\right)^{1/2}\\
&\lesssim&\left(\sum_{s=-\infty}^{\infty}\left(
\sum_{\bar s=-\infty}^\infty \sup_{z}
\sup_{\bar{s}\leq w<\bar{s}+1}H_{k,j}(|z|,|w|)
\int_{\dbR}\int_{s-\bar{s}}^{s-\bar{s}+1}\int_{-1}^1|g(z,y,\tau)|
d\tau\,dy\, dz\right)^{2}\right)^{1/2}.
\end{eqnarray*}
Using Minkowski inequality with respect to the sum on $\bar{s}$, we
continue with
\begin{eqnarray*}
&\lesssim&\sum_{\bar s=-\infty}^\infty \sup_{z}
\sup_{\bar{s}\leq w<\bar{s}+1}H_{k,j}(|z|,|w|)
\left(\sum_{s=-\infty}^{\infty}\left(
\int_{\dbR}\int_{s-\bar{s}}^{s-\bar{s}+1}\int_{-1}^1|g(z,y,\tau)|
d\tau\,dy\, dz\right)^{2}\right)^{1/2}\\
&\lesssim&\sum_{\bar s=-\infty}^\infty \sup_{z}
\sup_{\bar{s}\leq w<\bar{s}+1}H_{k,j}(|z|,|w|)
\left(\sum_{\sigma=-\infty}^{\infty}\left(
\int_{\dbR}\int_{\sigma}^{\sigma+1}\int_{-1}^1|g(z,y,\tau)|
d\tau\,dy\, dz\right)^{2}\right)^{1/2}\\
&\lesssim& 2^{2\alpha(k,j)}\left(\sum_{\sigma=-\infty}^{\infty}\left(
\int_{\dbR}\int_{\sigma}^{\sigma+1}\int_{-1}^1|g(z,y,\tau)|
d\tau\,dy\, dz\right)^{2}\right)^{1/2},
\end{eqnarray*}
where in the last step we used (\ref{maxx}). The proof is then
complete.
\end{proof}

\begin{remark}
If one wants to prove the classical $L^{2}_{x,y}$ estimate of the
maximal function, the price to pay is formally  an extra derivative
with respect to the $y$ variable.
\end{remark}

During the proof of the well-posedness result presented in Theorem
\ref{main1} we will need a {{\em weighted maximal function estimate}}.
The precise statement is presented in the following theorems

\begin{theorem}
   For any $\sigma>3/4, \gamma>1/2$ and $\theta>1,$ there exists
   $\alpha>1/2$ such that
\begin{eqnarray}
& &\label{ytmaxx}\left(\sum_{s=-\infty}^{\infty}\sup_{ s\leq y<s+1}
\sup_{|t|\leq 1}\sup_{x}\langle y\rangle^{\alpha}
|QU(t)u_0(x,y)|^2\right)^{1/2}\\
\nonumber&\lesssim& \|\langle y\rangle ^{\alpha}(1+D_x)^{\sigma}
(1+D_y)^{\gamma}u_0\|_{L^2(\dbR^2)}
+\|(1+D_x)^{\sigma-1/2}(1+D_y)^{\gamma+1/2}u_0\|_{L^2(\dbR^2)},\\
& &\label{ytmaxy}\left(\sum_{r=-\infty}^{\infty}\sup_{ r\leq x<r+1}
\sup_{|t|\leq 1}\sup_{y}\langle y\rangle^{\alpha}
|QU(t)u_0(x,y)|^2\right)^{1/2}\\
\nonumber&\lesssim& \|\langle y\rangle^{\alpha}(1+D_x)^{\sigma}
(1+D_y)^{\gamma}u_0\|_{L^2(\dbR^2)}+
\|\langle y\rangle^{\alpha}D_y^{\theta}u_0\|_{L^2(\dbR^2)}\\
\nonumber& &+\|(1+D_x)^{\sigma-1/2}D_y^{\gamma+1/2}
u_0\|_{L^2(\dbR^2)}
+\|(1+D_y)^{\theta+1/2}u_0\|_{L^2(\dbR^2)}.
\end{eqnarray}
\label{ymax}\end{theorem}
For low frequencies we have:
\begin{theorem}
For any $\sigma>1/4,  \gamma>1/2$ and $\theta>1$,
there exists $\alpha>1/2$ such that
\begin{eqnarray}
& &\label{ytmaxxlow}\left(\sum_{s=-\infty}^{\infty}\sup_{ s\leq y<s+1}
\sup_{|t|\leq 1}\sup_{x}\langle y\rangle^{\alpha}
|(Id-Q)U(t)u_0(x,y)|^2\right)^{1/2}\\
\nonumber&\lesssim& \|\langle y\rangle^{\alpha}D_x^{-\gamma}
(1+D_y)^{\gamma}u_0\|_{L^2(\dbR^2)}
+\|D_x^{-\theta}(1+D_y)^{\sigma}
u_0\|_{L^2(\dbR^2)}\\
& &\label{ytmaxylow}\left(\sum_{r=-\infty}^{\infty}\sup_{ r\leq x<r+1}
\sup_{|t|\leq 1}\sup_{y}\langle 
y\rangle^{\alpha}|(Id-Q)U(t)u_0(x,y)|^2\right)^{1/2}\\
\nonumber&\lesssim&
\|\langle y\rangle^{\alpha}D_x^{-\sigma}(1+D_y)^{\theta}
u_0\|_{L^2(\dbR^2)}
+\|D_x^{-\sigma-1/2}(1+D_y)^{\theta+1/2}u_0\|_{L^2(\dbR^2)}.
\end{eqnarray}
\label{ymaxlow}\end{theorem}
%
Here we prove only (\ref{ytmaxx}) of Theorem \ref{ymax}
and (\ref{ytmaxxlow}) of Theorem \ref{ymaxlow},
the rest follows  with similar arguments. 
\begin{proof}
Let $\tilde{Q}$ be the operator such that
$\F(\tilde{Q}f)(\xi,\lambda)=\chi_
{\{\left|\frac{\lambda}{\x}\right|\geq 1\}}
\hat{f}(\xi,\lambda)$. Assume $w_{0}(\x,\lambda)=w_{0}
\left(\left|\frac{\lambda}{\x}\right|\right)$,
where $w_{0}(r)$ is a smooth
characteristic function of the interval $[-2,2]$.  We write
$QU(t)u_{0}=(Id-\tilde{Q})QU(t)u_{0}+\tilde{Q}QU(t)u_{0}$.
Next we define the  multiplier operator $P_{k}$ such that
$\F(P_{k}f)(\x)=\psi(2^{-k}\xi)\hat{f}(\xi)$,  where
$\psi$ is a smooth characteristic function of $[1/2,2]$, and
$\sum_{k\geq 0}P_{k}=Q$.
We write $(Id-\tilde{Q})QU(t)u_{0}=\sum_{k\geq 0}(Id-\tilde{Q})U(t)
\sum P_{k}u_{0}$. We then define the operator
$$W_{k}^{z}=|y|^{z}(Id-\tilde{Q})U(t)P_{k}$$
for $z\in \dbC$ and we use complex interpolation. Assume
$\Re \,z=0$,
we use the argument presented to prove (\ref{tmaxx}) and we write
\begin{equation}
   \left(\sum_{s=-\infty}^{\infty}\sup_{ s\leq y<s+1}
\sup_{|t|\leq 1}\sup_{x}|(Id-\tilde{Q})U(t)P_{k}u_0(x,y)|^2
\right)^{1/2}\lesssim \|u_{0}\|_{L^{2}(2^{\alpha(k,j)}dx dy)},
\label{tildesigma=0}\end{equation}
where $\alpha(k,j)=5k/2+j+$.
Assume now that $\Re \,z=1$. We have
\begin{eqnarray*}
\F(y(Id-\tilde{Q})U(t)P_{k}u_0)(\x,\lambda)&\sim& \partial_{\lambda}
(e^{it(\x^{3}+\frac{\lambda^{2}}{\x})}
w_{0}(\lambda/\x)\psi(\x)\widehat{u_{0}})(\x,\lambda)\\
&\sim&
e^{it(\x^{3}+\frac{\lambda^{2}}{\x})}
\left(2\frac{\lambda}{\x}+\partial_{\lambda}\right)
(w_{0}(\lambda/\x)\psi(\x)\widehat{u_{0}})\\
&\sim& e^{it(\x^{3}+\frac{\lambda^{2}}{\x})}
\left(w_{0}(\lambda/\x)\psi(\x)\frac{\lambda}{\x}\widehat{u_{0}}
+2^{-k}\tilde{w}_{0}(\lambda/\x)\psi(\x)\widehat{u_{0}}+
w_{0}(\lambda/\x)\psi(\x)\partial_{\lambda}\widehat{u_{0}}\right)
\end{eqnarray*}
and agin using the argument presented to prove (\ref{tmaxx}) we have
\begin{equation}
   \left(\sum_{s=-\infty}^{\infty}\sup_{ s\leq y<s+1}
\sup_{|t|\leq 1}\sup_{x}|y(Id-\tilde{Q})U(t)P_{k}u_0(x,y)|^2
\right)^{1/2}\lesssim \|u_{0}\|_{L^{2}
(2^{\alpha(k,j)}(1+y^{2})dx dy)}.
\label{tildesigma=1}\end{equation}
We then  interpolate between (\ref{tildesigma=0})
and (\ref{tildesigma=1}) and sum over $k$ to obtain
(\ref{ytmaxx}) in this case.

Next define  $T_{k,j}, k\in \dbN, j \in \dbN$, such that
$\F(T_{k,j}(f))(\x,\lambda)=\theta_{k,j}\hat{f}(\xi,\lambda)$, where
$\theta_{k,j}$ was defined in Lemma \ref{I}, and
$\sum_{k,j}T_{k,j}^{2}=Id$. We first prove, again by interpolation,
that for any $\alpha\in [0,1]$
\begin{equation}
   \left(\sum_{s=-\infty}^{\infty}\sup_{ s\leq y<s+1}
\sup_{|t|\leq 1}\sup_{x}|y|^{\alpha}|\tilde{Q}QU(t)T_{k,j}u_0(x,y)|^2
\right)^{1/2}\lesssim \|u_{0}\|_{L^{2}(2^{\alpha(k,j)}(2^{2j}+
y^{2})^{\alpha}dx dy)}.
\label{dyadicymax}\end{equation}
where $\alpha(k,j)=5k/2+j+$. Define the operator
$$W_{k,j}^{z}=|y|^{z}\tilde{Q}QU(t)T_{k,j}.$$
Then  for $\Re \,z=0$,
the argument presented to prove (\ref{tmaxx}) gives
\begin{equation}
   \left(\sum_{s=-\infty}^{\infty}\sup_{ s\leq y<s+1}
\sup_{|t|\leq 1}\sup_{x}|\tilde{Q}QU(t)T_{k,j}u_0(x,y)|^2
\right)^{1/2}\lesssim \|u_{0}\|_{L^{2}(2^{\alpha(k,j)}dx dy)}.
\label{sigma=0}\end{equation}
For $\Re \, z=1$ we have
\begin{eqnarray}
\nonumber
\F(y\tilde{Q}QU(t)T_{k,j}u_0)(\x,\lambda)&\sim& \partial_{\lambda}
(e^{it(\x^{3}+\frac{\lambda^{2}}{\x})}\hat{u_{0}})(\x,\lambda)\\
\nonumber &\sim&
e^{it(\x^{3}+\frac{\lambda^{2}}{\x})}
\left(\frac{\lambda^{2}}{\x}+\partial_{\lambda}\right)\hat{u_{0}})
(\x,\lambda)\\
\label{2j}&\sim& 2^{j}\F(\tilde{Q}QU(t)T_{k,j}u_0)(\x,\lambda)+
\F(\tilde{Q}QU(t)T_{k,j}yu_0)(\x,\lambda)
\end{eqnarray}
Again, the argument presented to prove (\ref{tmaxx}) yields
\begin{equation}
   \left(\sum_{s=-\infty}^{\infty}\sup_{ s\leq y<s+1}
\sup_{|t|\leq 1}\sup_{x}|y\tilde{Q}QU(t)T_{k,j}u_0(x,y)|^2
\right)^{1/2}\lesssim \|u_{0}\|_{L^{2}
(2^{\alpha(k,j)}(2^{2j}+y^{2})dx dy)}.
\label{sigma=1}\end{equation}
Combining (\ref{sigma=0}) and (\ref{sigma=1}) and using complex
interpolation
we obtain (\ref{dyadicymax}). Now using (\ref{dyadicymax}) we have
\begin{eqnarray*}
& &\left(\sum_{s=-\infty}^{\infty}\sup_{ s\leq y<s+1}
\sup_{|t|\leq 1}\sup_{x}|y|^{\sigma}|\tilde{Q}QU(t)\sum_{k,j}T_{k,j}^{2}
u_0(x,y)|^2\right)^{1/2}\\&\lesssim& \sum_{k,j}\left(\sum_{s=-\infty}
^{\infty}\sup_{ s\leq y<s+1}
\sup_{|t|\leq 1}\sup_{x}|y|^{\sigma}|\tilde{Q}QU(t)T_{k,j}^{2}
u_0(x,y)|^2\right)^{1/2}\\
&\lesssim&\sum_{k,j}\|T_{k,j}u_{0}\|_{L^{2}(2^{\alpha(k,j)}(2^{2j}+
y^{2})^{\alpha}dx dy)}\\
&\lesssim&\sum_{k,j}(\|T_{k,j}u_{0}\|_{L^{2}(2^{\alpha(k,j)+
2\alpha j})
dx dy}+\|T_{k,j}u_{0}\|_{L^{2}(2^{\alpha(k,j)}|y|^{2\alpha}dx dy)}).
\end{eqnarray*}
For $\alpha=1/2+$ we have $\frac{\alpha(k,j)}{2}+\alpha j=
\frac{(1-)}{4}k+(1+)j$ and
\begin{equation}
\sum_{k,j}\|T_{k,j}u_{0}\|_{L^{2}(2^{\alpha(k,j)+2\sigma j}
dx dy}\lesssim \|D^{1/4}_{x}D^{1+}_{y}u_{0}\|_L^{2}.
\label{purel2}\end{equation}
To treat the term involving $|y|^{\alpha}$ we need a commutator
estimate. Recall  that here $j\geq 0$. Then for $\alpha>1/2,
\sigma=3/4+, \gamma=1/2+$ we have
$$
\|T_{k,j}u_{0}\|_{L^{2}(2^{\alpha(k,j)}|y|^{2\alpha}dx dy)}\sim
2^{\epsilon(k+j)}\||y|^{\alpha}T_{k,j}(D_{x}^{3/4+}
D_{y}^{1/2+}u_{0})\|_{L^{2}}.
$$
Write $g(x,y)=D_{x}^{3/4+}D_{y}^{1/2+}u_{0}$, then by
(\ref{linfinity}) it follows
\begin{eqnarray*}
\||y|^{\alpha}T_{k,j}g\|_{L^{2}}&\sim &\|D^{\alpha}_{\lambda}
\theta_{k,j}\hat{g}\|_{L^{2}}\\
&\lesssim& \|D^{\alpha}_{\lambda}(\theta_{k,j}\hat{g})
-D^{\alpha}_{\lambda}(\theta_{k,j})\hat{g}-
\theta_{k,j}D^{\alpha}_{\lambda}(\hat{g})\|_{L^{2}}\\
&+&\|D^{\alpha}_{\lambda}(\theta_{k,j})\|_{L^{2}}
\|\hat{g}\|_{L^{\infty}}+\|\theta_{k,j}\|_{L^{\infty}}
\|D^{\alpha}_{\lambda}(\hat{g})\|_{L^{2}}\\
&\lesssim&(\|D^{\alpha}_{\lambda}(\theta_{k,j})\|_{L^{2}}
+\|\theta_{k,j}\|_{L^{\infty}})
\|D^{\alpha}_{\lambda}(\hat{g})\|_{L^{2}}
\end{eqnarray*}
where, in the last step, we used Sobolev's theorem. Now, it is not
hard to show that $\|D^{\alpha}_{\lambda}(\theta_{k,j})\|_{L^{2}}
\sim 2^{(1/2-\alpha)(j+k)}$. In this case $j,k\geq 0$, hence
$\|D^{\alpha}_{\lambda}(\theta_{k,j})\|_{L^{2}}\leq C$, but in the
proof of (\ref{ytmaxxlow}) this ``small'' factor must be considered.
From here we obtain
\begin{equation}
\sum_{k,j} \|T_{k,j}u_{0}\|_{L^{2}(2^{\alpha(k,j)}
|y|^{1+}dx dy)})\lesssim \|D_{x}^{3/4+}D_{y}^{1/2+}u_{0}\|_{L^{2}}
+\|<y>^{1/2+}D_{x}^{3/4+}D_{y}^{1/2+}u_{0}\|_{L^{2}}.
   \label{yl2}\end{equation}
Then combining (\ref{purel2}) and (\ref{yl2}) we obtain
(\ref{ytmaxx}).

We now pass to the proof of (\ref{ytmaxxlow}). A weaker version of
this can be simply obtained following the argument above by
substituting \eqref{sigma=1} with
\begin{equation}
   \left(\sum_{s=-\infty}^{\infty}\sup_{ s\leq y<s+1}
\sup_{|t|\leq 1}\sup_{x}|y\tilde{Q}(Id-Q)U(t)T_{k,j}u_0(x,y)|^2
\right)^{1/2}\lesssim \|u_{0}\|_{L^{2}
(2^{\alpha(k,j)}(2^{2j}+y^{2})dx dy)},
\label{sigma=1low}\end{equation}
where $\alpha=j+$. But one can do better by repeating the proof of
\eqref{tmaxxlow} directly with $2^{j}\tilde{Q}QU(t)T_{k,j}u_0$,
that is the first contribution in \eqref{2j}. In this case,
\eqref{missingref} changes into
\begin{equation}
\tilde\delta(\tilde{k},\tilde{j})=
\max{(2\tilde{k}+\tilde{2j},-\chi_{\{\tilde{j}\geq 0\}}\tilde{j},
\chi_{\{\tilde{j}\geq 0\}}(5/2\tilde{k}+\tilde{2j}))},
\end{equation}
which translate into $\tilde\alpha(k,j)=j+$, if $2k+2j<j$ and
$\tilde\alpha(k,j)=2k+2j$, if $2k+2j\geq j$. The result then
follows by interpolation.

\end{proof}


\section{The Smoothing Effect Estimates}
In this section we prove some basic estimates that describe the
smoothing effect associated to the IVP (\ref{livp}). As mentioned in
the introduction, the smoothing effect is linked to the curvature of
the surface
$S=\left\{(\x,\lambda,\x^{3}+\frac{\lambda^{2}}{\x})\right\}$. We
show that in the region where the surface is of quadratic type
(i.e. $|\lambda|/|\x|\gg |\x|$), the dispersion behaves like in the
Schro\"odinger equation and we do not gain more that $\partial_{x}^{1/2}$
smoothness. On
the other had in the region where the surface is of cubic type (i.e.
$|\lambda|/|\x|\ll |\x|$), the dispersion behaves like in the KdV equation
and we gain the full derivative $\partial_{x}$.

\begin{lemma}
If $U(t)u_0$ is the solution of the IVP (\ref{livp}), then

\begin{equation}
\|\partial_x U(t)P_+Qu_0\|_{L^\infty_xL^2_{t,y}}\lesssim
\|u_0\|_{L^2_{x,y}},
\label{p+}\end{equation}
\begin{equation}
\|\partial_x^{1/2} U(t)P_-Qu_0\|_{L^\infty_yL^2_{t,x}}
\lesssim \|u_0\|_{L^2_{x,y}}.
\label{p-}\end{equation}
The inhomogeneous version of (\ref{p-}) becomes
\begin{equation}
\|\partial_x \int_0^tU(t-t')P_-Qf(x,y,t')\|_{L^\infty_yL^2_{t,x}}
\leq c \|f\|_{L^1_yL^2_{x,t}}.
\label{nhp-}\end{equation}
\end{lemma}

\begin{proof} Our proof follows the proof of
the one-dimensional  KdV smoothing effect presented in \cite{KPV2}.
To prove (\ref{p+}), we write
\begin{eqnarray*}
\partial_xU(t)P_{+}Qu_0(x,y)&=& \int_{\{|\xi|\geq
\max{(1,2|\lambda|/|\xi|)}\}}
e^{i(\xi,\lambda)(x,y)}\xi e^{it\phi(\xi,\lambda)}\widehat{u_0}
(\xi,\lambda)\,d\xi \,d\lambda\\
&=&\int \int_{\{|\xi|\geq
\max{(1,2|\lambda|/|\xi|)}\}}
e^{i\lambda y+it(\xi^3+\lambda^2/\xi)}[e^{i\xi x}\xi
\widehat{u_0}(\xi,\lambda)]\,d\xi\, d\lambda.
\end{eqnarray*}
We make the change of variables
$(\zeta, \lambda) = (\xi^3+\lambda^2/\xi,
\lambda)$, with  Jacobian $J(\xi,\lambda)\gtrsim |\xi|^2>1$.
We continue the chain of inequalities above with
\[
\int e^{i\lambda y+it\zeta}[e^{i\theta(\zeta,\lambda)x}
\theta(\zeta,\lambda)
\widehat{u_0}(\theta(\zeta,\lambda),\lambda)
\chi_{\{|\theta|\geq2|\lambda|/
|\theta|\}}J^{-1}]\,d\zeta d\lambda.\]
Then, by Plancherel
\begin{eqnarray*}
\|\partial_xU(t)P_{+}Qu_0\|_{L^2_{y,t}}&=&
\|e^{i\theta(\zeta,\lambda)x}
\theta(\zeta,\lambda)\widehat{u_0}(\theta(\zeta,\lambda),\lambda)
\chi_{\{|\theta|\geq2|\lambda|/|\theta|\}}J^{-1}
\|_{L^2_{\zeta,\lambda}}\\
&=&\left(\int|\theta(\zeta,\lambda)|^2|\widehat{u_0}|^2J^{-2}
\,d\zeta d\lambda
\right)^{1/2}\\
&\leq&\left(\int|\xi|^2|\widehat{u_0}|^2|\xi|^{-2}\right)^{1/2}
=\|u_0\|_{L^2_{x,y}}.
\end{eqnarray*}
To prove (\ref{p-}), we use a similar argument

\begin{eqnarray*}
\partial_x^{1/2}U(t)P_{+}Qu_0(x,y)&=&
\int_{\{1\leq |\xi|<2|\lambda|/|\xi|\}}
e^{i(\xi,\lambda)(x,y)}\xi^{1/2} e^{it\phi(\xi,\lambda)}
\widehat{u_0}(\xi,\lambda)\,d\xi d\lambda\\
&=&\int e^{i\xi x+it(\xi^3+\lambda^2/\xi)}[e^{i\lambda y}\xi^{1/2}
\widehat{u_0}(\xi,\lambda)]\,d\xi d\lambda.
\end{eqnarray*}
We make the change of variables
$(\xi,\rho) = (\xi,\xi^3+\lambda^2/\xi)$,
with Jacobian $J(\xi,\lambda)\gtrsim |\xi|$.
We set $\lambda=\gamma(\xi,\rho)$ and we continue the chain of
inequalities above with
\[
\int e^{i\xi x+it\rho}[e^{i\gamma(\zeta,\lambda)y}\xi^{1/2}
\widehat{u_0}(\xi,\gamma(\zeta,\lambda))\chi_{\{|\xi|<2|\gamma|/
|\theta|\}}J^{-1}]\,d\rho d\xi.\]
Then, by Plancherel
\begin{eqnarray*}
\|\partial_x^{1/2}U(t)P_{-}Qu_0\|_{L^2_{x,t}}&=&
\|e^{i\gamma(\zeta,\lambda)y}
\xi^{1/2}\widehat{u_0}(\xi,\gamma(\zeta,\lambda))
\chi_{\{|\xi|<2|\rho|/|\theta|\}}J^{-1}\|_{L^2_{\xi,\rho}}\\
&=&\left(\int|\xi||\widehat{u_0}|^2J^{-2}\,d\xi d\rho
\right)^{1/2}\\
&\leq&\left(\int|\xi||\widehat{u_0}|^2|\xi|^{-1}\right)^{1/2}
=\|u_0\|_{L^2_{x,y}}.
\end{eqnarray*}
In order to prove the inhomogeneous estimate (\ref{nhp-}),
we first observe that the
adjoint operator $(U(t)P_{-}Q)^{*}:L^1_{y}L^2_{x,t}\longrightarrow
L^2_{x,y}$ is defined as

\[
(U(t)P_{-}Q)^{*}f(x,y,t)=\int_{-\infty}^\infty P_-U(-t)Qf(x,y,t)\,dt
\]
and that  the dual version of (\ref{nhp-}) becomes

\begin{equation}
\|\partial_x^{1/2}\int_{-\infty}^\infty P_-U(-t)Qf(x,y,t)\,dt
\|_{L^2_{x,y}}
\leq \|f\|_{L^1_yL^2_{x,t}}.
\end{equation}
Now, following the argument on page 554 of \cite{KPV2},
we see that it is enough to show

\begin{equation}
\|\partial_x\int_{-\infty}^\infty U(t-t')P_-fQ(x,y,t')\,dt'\|
_{L^\infty_yL^2_{x,t}}\leq \|f\|_{L^1_yL^2_{x,t}}.
\end{equation}
To see this we take a smooth function $g$ such that
$\|g\|_{L^1_yL^2_{x,t}}\leq 1$ and write

\begin{eqnarray*}
& &\int_{-\infty}^{\infty}\int_{\dbR^2} 
\left(\partial_x\int_{-\infty}^\infty U(t-t')
P_-fQ(x,y,t')\,
dt'\right)g(x,y,t)\,dx\, dy\, dt\\
&=&\int_{-\infty}^\infty \int_{-\infty}^\infty \int_{\dbR^2}
\xi e^{i(t-t')\phi(\xi,\lambda)}
\chi_{\{1\leq |\xi|<2|\lambda|/|\xi|\}}\hat{f}(\xi,\lambda,t')
\check{g}(\xi,\lambda,t')\,d\xi\,d\lambda \,dt\,dt'\\
&=&\int_{-\infty}^\infty \int_{-\infty}^\infty \int_{\dbR^2}
\xi^{1/2} e^{-it'\phi(\xi,\lambda)}
\chi_{\{1\leq |\xi|<2|\lambda|/|\xi|\}}\hat{f}(\xi,\lambda,t')
\xi^{1/2} e^{it\phi(\xi,\lambda)}\check{g}(\xi,\lambda,t)
\,d\xi\,d\lambda \,dt\,dt'\\
&=&\int_{\dbR^2}\partial_x^{1/2}\int_{-\infty}^\infty
U(-t')P_-Qf(x,y,t')\,dt'\partial_x^{1/2}\int_{-\infty}^\infty
U(-t)P_-Qg(x,y,t)\,dt\,dx\, dy\\
&\leq&\left\|\partial_x^{1/2}\int_{-\infty}^\infty
U(-t')P_-Qf(x,y,t')\,dt'\right\|_{L^2_{x,y}}\left\|\partial_x^{1/2}
\int_{-\infty}^\infty U(-t)P_-Qg(x,y,t)\,dt\right\|_{L^2_{x,y}}\\
&\leq&\|f\|_{L^1_yL^2{x,y}}\|g\|_{L^1_yL^2{x,y}}.
\end{eqnarray*}
\end{proof}


\label{group}\section{The Group Estimates}

In order to define a Banach space suitable for  a
contraction mapping theorem,  we need to
analyze how the group $U(t)$ acts in certain weighted anisotropic Sobolev
spaces.

\begin{lemma}
Assume that $<y>^\alpha g, \, (1+D_x)^{-\alpha}(1+D_y)^\alpha g
\in L^2_{x,y}$ for $\alpha \in [0,1]$. Then

\begin{equation}
\|\langle y\rangle^\alpha U(t)Qg(x,y)\|_{L^2_{x,y}}\leq
\|\langle y\rangle^\alpha g\|_{L^2_{x,y}}+
\|(1+D_y)^\alpha (1+D_x)^{-\alpha})g\|_{L^2_{x,y}}.
\label{ygroup}\end{equation}
For the small frequencies we also have
\begin{equation}
\|\langle y\rangle^\alpha U(t)(Id-Q)g(x,y)\|_{L^2_{x,y}}\leq
\|\langle y\rangle^\alpha g\|_{L^2_{x,y}}+
\|(1+D_y)^\alpha D_x^{-\alpha-\epsilon})g\|_{L^2_{x,y}}.
\label{ygrouplow}\end{equation}
\end{lemma}
\begin{proof}
Here we only present the proof for (\ref{ygroup}), the one for
(\ref{ygrouplow}) follows by similar arguments.
We define
$\tilde{Q}$ to be the operator such that
$\F(\tilde{Q}f)(\lambda)=\chi_
{\{|\lambda|\geq 1\}}\hat{f}(\lambda)$. (Note that this $\tilde{Q}$
is different from the one used in the proof of Theorem 3.7.) 
Then assume
$\psi_{0}(\lambda)=\psi_{0}(|\lambda|)$, is a smooth
characteristic function of the interval $[-2,2]$. We consider first
the operator $T^{z}=\langle y\rangle^{z}U(t)(Id-\tilde{Q})
Qg(x,y)$, for $z \in \dbC$. Clearly for $\Re \,z=0$
\begin{equation}
   \|U(t)(Id-\tilde{Q})Qg(x,y)\|_{L^{2}}\lesssim
   \|g\|_{L^{2}(dx dy)}.
\label{lows=0}\end{equation}
On the other hand, when $\Re \, z=1$ we estimate
\begin{eqnarray*}
\F(yU(t)(Id-\tilde{Q})Qg)(\x,\lambda)&\sim& \partial_{\lambda}
(e^{it(\x^{3}+\frac{\lambda^{2}}{\x})}\psi_{0}(\lambda)\hat{g})\\
&\sim&e^{it(\x^{3}+\frac{\lambda^{2}}{\x})}\left(\frac{\lambda}{\x}
\hat{g}+w'(\lambda)\hat{g}+\partial_{\lambda}
\hat{g}\right)\chi_{\{|\x|\geq 1\}}.
\end{eqnarray*}
It follows that for $\Re \, z=1$
\begin{equation}
   \sup_{|t|\leq 1}\|yU(t)(Id-\tilde{Q})Qg(x,y)\|_{L^{2}}\lesssim
   \|g\|_{L^{2}((1+y^{2})dx dy)}.
\label{lows=1}\end{equation}
Interpolating (\ref{lows=0}) and (\ref{lows=1}) we obtain
$$
\sup_{|t|\leq 1}\|\langle y\rangle^{\alpha}U(t)(Id-\tilde{Q})Qg(x,y)
\|_{L^{2}}\lesssim
   \|g\|_{L^{2}((1+y^{2})^{\alpha}dx dy)},
$$
for $\alpha\in [0,1]$,
and this proves (\ref{ygroup}) for $|\lambda|\leq 1$.
To treat the region $|\lambda|> 1$
we introduce a Littlewood-Paley decomposition. Let
$\psi_{k,j}, k,j\in \dbN$, such that
$$\psi_{k,j}(\x,\lambda)=\psi(2^{-k}\x,2^{-j}\lambda)$$
where $\psi(|r|,|s|)$ is a smooth characteristic function of the
rectangle $[1/2,2]\times[1/2,2]$. Let
$P_{k,j}$ be the multiplier operator such that
\begin{equation}\F(P_{k,j}(f))(\x,\lambda)=\psi_{k,j}
   \hat{f}(\x,\lambda).
\label{pkj}\end{equation}
We need
some preliminary estimates:
\begin{equation}
\left\||y|^{\alpha}(\sum_{k,j}|P_{k,j}f|^{2})^{1/2}
\right\|_{L^{2}}
\lesssim \||y|^{\alpha}f\|_{L^{2}} + \|f\|_{L^{2}},
\label{y1}\end{equation}
\begin{equation}
\left \||y|^{\alpha}(\sum_{k,j}|P_{k,j}f_{k,j}|)\right\|_{L^{2}}
\lesssim \left\||y|^{\alpha}(\sum_{k,j}|f_{k,j}|^{2})^{1/2}
\right\|_{L^{2}},
\label{y2}\end{equation}
\begin{equation}
\left\|\frac{1}{\langle y\rangle^{\alpha}}(\sum_{k,j}|P_{k,j}f|^{2})^{1/2}
\right\|_{L^{2}}
\lesssim \left\|\frac{1}{\langle y\rangle^{\alpha}}f\right\|_{L^{2}},
\label{y3}\end{equation}
\begin{equation}
\|\langle y\rangle^{\alpha}f\|_{L^{2}}\lesssim \left
\|\langle y\rangle^{\alpha}(\sum_{k,j}
|P_{k,j}f|^{2})^{1/2}\right\|_{L^{2}}.
\label{y4}\end{equation}
To prove (\ref{y1}), note that the square of the L.H.S. is 
$\sum_{k,j}\||y|^{\alpha}P_{k,j}f\|_{L^{2}}^{2}$. We
estimate this by interpolation. Define the operator
$T_{k,j}^{z}=|y|^{z}P_{k,j}$. Then,
when $\Re \,z=0$ we get the trivial
bound $\|P_{j,k}f\|_{L^{2}}^{2}$. When $\Re z=1$ we get
\begin{eqnarray*}
\|yP_{k,j}f\|^{2}_{L^{2}}&=&\|\partial_{\lambda}
(\psi_{k,j}f)\|_{L^{2}}^{2}\\
&\lesssim &2^{-2j}\|\psi'_{k,j}f\|_{L^{2}}^{2}+
\|\psi_{k,j}\partial_{\lambda}f\|_{L^{2}}^{2}\\
&\lesssim&2^{-2j}\|\tilde{P}_{k,j}f\|_{L^{2}}^{2}+
\|P_{k,j}(yf)\|_{L^{2}}^{2}.
\end{eqnarray*}
We now use the support properties of $\tilde{P}_{k,j}$ and $P_{k,j}$
on the Fourier transform to conclude that
when $\Re \,z=0$
\begin{equation}
   \left\||y|^{z}\left(\sum_{k,j}|P_{k,j}f|^{2}\right)^{1/2}
   \right\|_{L^{2}}\lesssim \|f\|_{L^{2}},
\label{y1s=0}\end{equation}
and when $\Re \,z=1$
\begin{equation}
   \left\||y|^{z}\left(\sum_{k,j}|P_{k,j}f|^{2}\right)^{1/2}
   \right\|_{L^{2}}\lesssim \|\langle y\rangle f\|_{L^{2}}.
\label{y1s=1}\end{equation}
Complex interpolation between (\ref{y1s=0}) and (\ref{y1s=1}) proves
(\ref{y1}). We now prove (\ref{y2}). We use again interpolation.
When $\Re \,z=0$:
\begin{eqnarray*}
\left\|\sum_{k,j}P_{k,j}f_{k,j}\right\|_{L^{2}}&=&
\sup_{\|g\|_{L^{2}}\leq 1}
\left|\sum_{k,j}\int P_{k,j}f_{k,j}g \, dx\, dy\right|
=\sup_{\|g\|_{L^{2}}\leq 1}\left|\sum_{k,j}
\int f_{k,j}P_{k,j}g \,dx\, dy\right|\\
&\leq& \sup_{\|g\|_{L^{2}}\leq 1}\int
\left(\sum_{k,j}|f_{k,j}|^{2}\right)^{1/2}
\left(\sum_{k,j}|P_{k,j}g|^{2}\right)^{1/2}\, dx\, dy
\lesssim\left\|\left(\sum_{k,j}|f_{k,j}|^{2}
\right)^{1/2}\right\|_{L^{2}}
   \end{eqnarray*}
where, in the last step, we used Littlewood-Paley theory.
When $\Re \,z=1$ we use similar arguments to show that
\begin{eqnarray*}
\left\|y\sum_{k,j}P_{k,j}f_{k,j}\right\|_{L^{2}}&=&
\left\|\sum_{k,j}2^{-j}\tilde{P}_{k,j}f_{k,j}+
\sum_{k,j}P_{k,j}(yf_{k,j})\right\|_{L^{2}}\\
&\lesssim& \left\|\left(\sum_{k,j}|f_{k,j}|^{2}\right)^{1/2}
\right\|_{L^{2}} + \left\|\left(\sum_{k,j}|yf_{k,j}|^{2}\right)^{1/2}
\right\|_{L^{2}}\\
&\lesssim& \left\|\langle y\rangle\left(\sum_{k,j}|f_{k,j}|^{2}\right)^{1/2}
\right\|_{L^{2}}.
\end{eqnarray*}
We are now ready to prove (\ref{y3}). Let $g_{jk}$ be  functions
such that
$$
\left\|\left(\sum_{k,j}|g_{k,j}|^{2}\right)^{1/2}\langle
y\rangle^{\alpha}\right\|_{L^{2}}\leq 1.
$$
Then
\begin{eqnarray*}
   \left\|\frac{1}{\langle y\rangle^{\alpha}}(\sum_{k,j}
   |P_{k,j}f|^{2})^{1/2}
\right\|_{L^{2}}&=& \sup_{g_{k,j}}
\left|\int\sum_{k,j}P_{k,j}fg_{k,j}\,dxdy\right|\\
&=& \sup_{g_{k,j}}\left|\int\sum_{k,j}P_{k,j}g_{k,j}f\,dxdy\right|\\
&\lesssim& \left\|\langle y\rangle^{\alpha}\sum_{k,j}P_{k,j}g_{k,j}
\right\|_{L^{2}}\left\|\frac{1}{\langle y\rangle^{\alpha}}
f\right\|_{L^{2}}.
\end{eqnarray*}
We then use (\ref{y2}) to obtain the desired result. We now prove
(\ref{y4}).  Write $f=\sum_{k,j}P^{2}_{k,j}f$. Let $g$ such that
$\|(1+|y|)^{-\alpha}g\|_{L^{2}}\leq 1$. We estimate
\begin{eqnarray*}
\int f(x,y)g(x,y)\, dx\, dy&=& \int \sum_{k,j}P_{k,j}^{2}f
g\,dx\,dy\\
&=& \int\sum_{k,j}P_{k,j}f P_{k,j}g\, dx\, dy\leq
\int\left(\sum_{k,j}|P_{k,j}f|^{2}\right)^{1/2}
\left(\sum_{k,j}|P_{k,j}g|^{2}\right)^{1/2}\, dx\, dy\\
&\lesssim& \left\|\langle y\rangle^{\alpha}\left(\sum_{k,j}
|P_{k,j}f|^{2}\right)^{1/2}\right\|_{L^{2}}
\left\|\langle y\rangle^{-\alpha}\left(\sum_{k,j}
|P_{k,j}g|^{2}\right)^{1/2}\right\|_{L^{2}}\\
&\lesssim&\left\|\langle y\rangle^{\alpha}\left(\sum_{k,j}
|P_{k,j}f|^{2}\right)^{1/2}\right\|_{L^{2}}
\|\langle y\rangle^{-\alpha}g\|_{L^{2}},
 \end{eqnarray*}
where we used (\ref{y3}) in the last step.
We are now ready to prove (\ref{ygroup}). Using the same ideas we
presented above it is not hard to show that
for any $\alpha \in [0,1]$
\begin{equation}
\sup_{|t|\leq 1}\|\langle y\rangle^{\alpha}P_{k,j}
\tilde{Q}QU(t)u_{0}\|_{L^{2}}
\lesssim \|u_{0}\|_{L^{2}((2^{2j-2k}+|y|^{2})^{\alpha})dx dy)}.
\label{dyadicygroup}\end{equation}
Then, by (\ref{y4}), we write
\begin{eqnarray}
\nonumber\sup_{|t|\leq 1}\|\langle y\rangle^{\alpha}\tilde{Q}
QU(t)u_{0}\|_{L^{2}}^{2}
&\lesssim&\|\langle y\rangle^{\alpha}(\sum_{0\leq 
k,j}|P_{k,j}U(t)u_{0}|^{2})^{1/2}
\|_{L^{2}}^{2}\\
&=&\sum_{0\leq k,j}\|\langle y\rangle^{\alpha}|P_{k,j}U(t)u_{0}|
\|_{L^{2}}^{2}
\label{step1}\end{eqnarray}
Let $\tilde{\psi}_{k,j}$ be such that
$\tilde{\psi}_{k,j}\psi_{k,j}=\psi_{k,j}$. We continue (\ref{step1}) with
$$
\sum_{0\leq k,j}\|\langle y\rangle^{\alpha}|U(t)P_{k,j}\tilde{P}_{k,j}u_{0}|
\|_{L^{2}}^{2}.
$$
Recall also that $\|U(t)P_{k,j}f\|^{2}_{L^{2}}\lesssim
\|f\|^{2}_{L^{2}}$. Thus, our basic estimate (\ref{dyadicygroup})
gives
\begin{eqnarray*}
\sup_{|t|\leq 1}\|\langle y\rangle^{\alpha}\tilde{Q}QU(t)
u_{0}\|_{L^{2}}^{2}
&\lesssim&\sum_{0\leq k,j}\|\langle 
y\rangle^{\alpha}|U(t)P_{k,j}\tilde{P}_{k,j}u_{0}|
\|_{L^{2}}^{2}\\
&\lesssim& \sum_{0\leq k,j}\|\tilde{P}_{k,j}u_{0}\|
_{L^{2}((2^{2j-2k}+|y|^{2})^{\alpha})dx dy)}^{2}\\
&\lesssim&\sum_{0\leq k,j}(\|2^{j\alpha-k\alpha}
\tilde{P}_{k,j}u_{0}\|
_{L^{2}}^{2}+
\||y|^{\alpha}\tilde{P}_{k,j}u_{0}\|_{L^{2}}^{2}\\
&\lesssim& \|(1+D_{x})^{-\alpha}(1+D)^{\alpha}_{y}u_{0}\|_{L^{2}}^{2}
+\||y|^{\alpha}(\sum_{0\leq k,j}(\tilde{P}_{k,j}u_{0})^{2})^{1/2}\|
_{L^{2}}^{2}\\
&\lesssim& \|(1+D_{x})^{-\alpha}(1+D_{y})^{\alpha}u_{0}
\|_{L^{2}}^{2}+\|\langle y\rangle^{\alpha}u_{0}\|_{L^{2}}^{2} ,
\end{eqnarray*}
where, in the last step, we used (\ref{y1}). This concludes the
proof of (\ref{ygroup}) and hence the proof of the lemma.
\end{proof}
We recall the Strichartz estimates associated to the IVP
(\ref{ivp}) due to Ben-Artzi and Saut \cite{BS}:
\begin{proposition}
If $0\leq \theta<1$, and $(q,p)=\left(\frac{2}{1-\theta},
\frac{2}{\theta}\right)$,  then
\begin{equation}
   \|U(t)u_{0}\|_{L^{q}_{t}L^{p}_{x,y}}\lesssim
   \|u_{0}\|_{L^{2}}
\label{xyinfinity}\end{equation}
\end{proposition}
In order to be able to estimate the intermediate orders of
derivatives in the non-linear term of (\ref{ivp}),
we need (\ref{xyinfinity}) with $(q,p)=(2,\infty)$
and with the weight $\langle y\rangle^{\alpha}, \alpha=1/2+$.
For this purpose we prove the following proposition:
\begin{proposition}
For any $\alpha \in [0,1]$ and $t\in [0,1]$
\begin{equation}
   \|\langle y\rangle^{\alpha}QU(t)u_{0}\|_{L^{2}_{t}
   L^{\infty}_{x,y}}\lesssim T^{\delta(\epsilon)}
   (\|\langle y\rangle^{\alpha}((1+D_{x})^{\epsilon}+
   (1+D_{y})^{\epsilon})
   Qu_{0}\|_{L^{2}}+\|(1+D_{y})^{\alpha+\epsilon}Qu_{0}\|_{L^{2}}).
\label{yxyinfinity}\end{equation}
For the low frequencies we have
\begin{equation}
   \|\langle y\rangle^{\alpha}(Id-Q)U(t)u_{0}\|_{L^{2}_{t}
   L^{\infty}_{x,y}}\lesssim
   T^{\delta(\epsilon)}(\|\langle y\rangle^{\alpha}
   (1+D_{y})^{\epsilon}
   u_{0}\|_{L^{2}}+\|D_{x}^{-\alpha-\epsilon}
   (1+D_{y})^{\alpha+2\epsilon}(Id-Q)u_{0}\|_{L^{2}})
\label{yxyinfinitylow}\end{equation}
for any $0<\epsilon<1$.
\end{proposition}
\begin{proof}
Here we prove only (\ref{yxyinfinity}); (\ref{yxyinfinitylow})
will follow from similar arguments. We first observe that
\begin{equation}
   \|\langle y\rangle^{\alpha}QU(t)u_{0}\|_{L^{2}_{t}
   L^{\infty}_{x,y}}\lesssim T^{\delta(\epsilon)}
\|[(1+D_{x})^{\epsilon}+(1+D_{y})^{\epsilon}]
\langle y\rangle^{\alpha}QU(t)u_{0}\|_{L^{q(\epsilon)}_{t}
   L^{p(\epsilon)}_{x,y}}
\label{sobolev}\end{equation}
where $q(\epsilon)=\frac{2}{\epsilon}, p(\epsilon)=
\frac{2}{1-\epsilon}, \delta(\epsilon)=\frac{\epsilon}{2}$.
Then by (\ref{weight}) we can continue with
$$T^{\delta(\epsilon)}
\|\langle y\rangle^{\alpha}
[(1+D_{x})^{\epsilon}+(1+D_{y})^{\epsilon}]QU(t)u_{0}\|_{L^{q(\epsilon)}_{t}
L^{p(\epsilon)}_{x,y}}
$$
and from here on we can proceed by interpolation, like in the proof
of Theorem \ref{ymax}, using (\ref{xyinfinity}).
\end{proof}

\label{fpa}\section{Proof of Theorems \ref{main1} and \ref{main2}:
The fixed point argument}
In this section we give only an outline of the proof of Theorem
\ref{main1}. The basic argument is illustrated in many articles
in the literature of dispersive equations. In particular, we refer
the reader to \cite{KPV3}. The proof of Theorem \ref{main2} is similar
and the main difference is explained in Remark \ref{lwp} below.

We start by transforming \eqref{ivp2} into the system of
integral equations
\begin{eqnarray}
\label{d1} u_{1}&=&\chi_{[0,T]}U(t)w_{0}+\chi_{[0,T]}\int_{0}^{t}
QU(t-t')\partial_{x}(u_{1}^{2}+u_{2}^{2}+2u_{1}u_{2})dt'\\
\label{d2}u_{2}&=&\chi_{[0,T]}U(t)v_{0}+\chi_{[0,T]}\int_{0}^{t}
(Id-Q)U(t-t')\partial_{x}(u_{1}^{2}+u_{2}^{2}+2u_{1}u_{2})dt'.
\end{eqnarray}
Then it is clear that a solution for \eqref{d1} and \eqref{d2}
is a fixed point for the operator
\begin{equation}
    \label{oper}L(w,v)=(L_{1}(w,v),L_{2}(w,v))
    \end{equation}
where
\begin{eqnarray}
 \label{l1}L_{1}(w,v)&=&\chi_{[0,T]}U(t)w_{0}+\chi_{[0,T]}\int_{0}^{t}
QU(t-t')\partial_{x}(w^{2}+v^{2}+2wv)dt'\\
\label{l2}L_{2}(w,v)&=&\chi_{[0,T]}U(t)v_{0}+\chi_{[0,T]}\int_{0}^{t}
(Id-Q)U(t-t')\partial_{x}(w^{2}+v^{2}+2wv)dt'.
\end{eqnarray}
We now define the  Banach space $X_{T}\times Y_{T}$ where
we will find the fixed point for (\ref{oper}).
\begin{definition}
Let $\sigma_{1}>3/4,
\gamma_{1}>1/2, \sigma_{2}>1,\gamma_{2}>1$ and let
$\alpha=\alpha(\sigma_{1},\sigma_{2},\gamma_{1},\gamma_{2})=1/2+$
to be the smallest $\alpha$ that can be chosen in Theorems
\ref{ymax} and \ref{ymaxlow}. Consider the  norms
\begin{eqnarray*}
   \|v\|_{1}&=& \|v\|_{H^{2}}\\
   \|v\|_{2}&=&\|\langle y\rangle^{\alpha}(1+D_{x})^{\sigma_{1}}
   (1+D_{y})^{\gamma_{1}}v\|
   _{L^{\infty}_{t}L^{2}_{x,y}}\\
\|v\|_{3}&=&\|\langle y\rangle^{\alpha}
   (1+D_{y})^{\gamma_{2}}v\|_{L^{\infty}_{t}L^{2}_{x,y}}\\
   \|v\|_{4}&=&\|\langle y\rangle^{\alpha}(1+D_{x})^{\sigma_{2}}u\|
   _{L^{\infty}_{t}L^{2}_{x,y}}\\
   \|v\|_{5}&=&\|\partial_{x}P_{+}(D_{x}^{2}
   +D_{y}^{2})v\|_{L^{\infty}_{x}L^{2}_{t,y}}\\
 \|v\|_{6}&=&\|\partial_{x}P_{+}(1+D_{x})^{\sigma_{1}}
   (1+D_{y})^{\gamma_{1}}v\|_{L^{\infty}_{x}L^{2}_{t,y}}\\
\|v\|_{7}&=&\|\partial_{x}P_{-}(D_{x}^{2}
   +D_{y}^{2})v\|_{L^{\infty}_{y}L^{2}_{t,x}}\\
\|v\|_{8}&=&\|\partial_{x}P_{-}(1+D_{x})^{\sigma_{1}}
   (1+D_{y})^{\gamma_{1}}v\|_{L^{\infty}_{y}L^{2}_{t,x}}\\
  \|v\|_{9}&=&\|\langle y\rangle^{\alpha}v\|_{L^{2}_{y}
L^{\infty}_{x,t}}\\
\|v\|_{10}&=&\|\langle y\rangle^{\alpha}v\|_{L^{2}_{x}
L^{\infty}_{y,t}}\\
\|v\|_{11}&=&\|\langle y\rangle^{\alpha}\partial_{y}v\|
   _{L^{2}_{t}L^{\infty}_{x,y}}\\
   \|v\|_{12}&=&\|\langle y\rangle^{\alpha}\partial_{x}v\|
   _{L^{2}_{t}L^{\infty}_{x,y}},
\end{eqnarray*}
and
\begin{equation}
   \|v\|_{X}=\max_{i=1\ldots 12}(\|v\|_{i}).
\label{normx}\end{equation}
Note that norms $1, \dots , 4$ relate to the group estimates, 
norms $5, \dots 8$ relate to the smoothing estimates and norms $9, 10, 11$
concern the maximal function estimates. 
Below we denote with $\|\cdot\|_{i,t}, i=1,\ldots, 4$ the norm
obtained from $\|\cdot\|_{i}, \, i=1,\ldots, 4$ when the time
variable is left free.
Recall that $F$ denotes the Fourier transform. 
We define the spaces of functions
$$X_{T}=
\{f(x,y,t) / t\in [0,T], \F_{x}(f)(\xi) = 0 \mbox{ if } |\x|\leq 1,
\mbox{ and } \|f\|_{X}<\infty\}$$
and
$$Y_{T}=\{f(x,y,t)/ t\in [0,T], \F_{x}(f)(\xi) = 0 \mbox{ if }
|\x|>1 \mbox{ and } \|f\|_{Y}<\infty\}$$
where
$$\|f\|_{Y}=\max_{i=1,3,9,10,11,12}\|f\|_{i}.$$
We  combine $X_{T}$ and $Y_{T}$ to obtain the space
$$Z_{T}=\{f(x,y,t)/ t\in [0,T], \|Qf\|_{X}+\|(Id-Q)f\|_{Y}<\infty\}.$$
\label{spacexy}\end{definition}
We also need  a space for the initial data. For this we introduce the
following definition:
\begin{definition}
For any $g(x,y)$ we define
$$\|g\|_{X,0}=\max_{i=1,\ldots , 4}\|v_{0}\|_{i},$$
(where we ignore the $L^\infty_{t}$ part of the $i= 1, \dots, 4$ norms since
we are considering initial data) 
and the corresponding space
$$X_{0}=
\{g(x,y) / \F_{x}(g)(\xi) = 0 \mbox{ if } |\x|\leq 1,
\mbox{ and } \|f\|_{X,0}<\infty\}.$$
We also define
$$\|g\|_{Y,0}=\max(
\|g\|_{L^{2}}, \, \,
\|\langle y\rangle^{\alpha}D_{x}^{-1/4-}
D_{y}^{1+}g\|_{L^{2}}, \, \,\|\langle y\rangle^{\alpha}
D_{x}^{-1/2-}
D_{y}^{1/2+}g\|_{L^{2}}, \, \, \|D_{x}^{-3/4-}
D_{y}^{3/2+}g\|_{L^{2}}<)\infty,$$
and the associated space
$$Y_{0}=
\{g(x,y) / \F_{x}(g)(\xi) = 0 \mbox{ if } |\x|> 1,
\mbox{ and } \|f\|_{Y,0}<\infty\}.$$
We  combine $X_{0}$ and $Y_{0}$ to obtain the space
$$Z_{0}=\{g(x,y)/ \|Qg\|_{X,0}+\|(Id-Q)g\|_{Y,0}<\infty\}.$$
\label{spacexy0}
\end{definition}
We now show that if
\begin{equation}
    \max(\|w_{0}\|_{X,0}, \|v_{0}\|_{X,0})=\ll \delta,
\label{smalldelta}\end{equation}
for an appropriate $\delta$, then the operator $L$ is a smooth
contraction in a ball of
$X_{T}\times Y_{T}$, centered at the origin of radius
$R\sim \delta$. First, we observe that, by definition,
$\F_{x}(L_{1}(w,v))(\xi)=0$ for $|\xi|\leq 1$ and
$\F_{x}(L_{2}(w,v))(\xi)=0$ for $|\xi|> 1$. 
We start with the estimate of the norm of $L_{1}(w,v)$.
We decided not to write explicitly the estimates for its linear
term because they are basically contained in Section 5.
We then start by  estimating the nonlinear term of \eqref{l1}
containing the  term $w^{2}$.
From  (\ref{ygroup}), \eqref{p+}, \eqref{nhp-},  
(\ref{ytmaxx}), (\ref{ytmaxy}) and
(\ref{yxyinfinity}) we have
\begin{eqnarray*}
   \|\int_{0}^{t}U(t-t')\partial_{x}w^{2}(t')dt'\|_{1,t}&\lesssim&
   \int_{0}^{T}\|\partial_{x}w^{2}(t')\|_{1,t'}dt'\\
\|\int_{0}^{t}U(t-t')\partial_{x}w^{2}(t')dt'\|_{2,t}&\lesssim&
   \int_{0}^{T}(\|\partial_{x}w^{2}(t')\|_{2,t'}+
\|\partial_{x}w^{2}(t')\|_{1,t'})dt'\\
\|\int_{0}^{t}U(t-t')\partial_{x}w^{2}(t')dt'\|_{3,t}&\lesssim&
   \int_{0}^{T}(\|\partial_{x}w^{2}(t')\|_{3,t'}+
\|\partial_{x}w^{2}(t')\|_{1,t'})dt'\\
\|\int_{0}^{t}U(t-t')\partial_{x}w^{2}(t')dt'\|_{4,t}&\lesssim&
   \int_{0}^{T}(\|\partial_{x}w^{2}(t')\|_{4,t'}+
\|\partial_{x}w^{2}(t')\|_{1,t'})dt'\\
\|\int_{0}^{t}U(t-t')\partial_{x}w^{2}(t')dt'\|_{5}&\lesssim&
   \int_{0}^{T}\|\partial_{x}w^{2}(t')\|_{1,t'}dt'\\
\|\int_{0}^{t}U(t-t')\partial_{x}w^{2}(t')dt'\|_{6}&\lesssim&
   \int_{0}^{T}\|\partial_{x}w^{2}(t')\|_{2,t'}dt'\\
\|\int_{0}^{t}U(t-t')\partial_{x}w^{2}(t')dt'\|_{7}&=&
\|\partial_{x}\int_{0}^{t}P_{-}U(t-t')(D_{x}^{2}
+D_{y}^{2})\partial_{x}(w^{2})(t')
dt'\|_{L^{\infty}_{y}L^{2}_{x,t}}\\
&\lesssim&\|\partial_{x}(D_{x}^{2}
+D_{y}^{2})(w^{2})\|_{L^{1}_{y}L^{2}_{x,t}}
\lesssim (\int\|\langle y\rangle^{\alpha}
\partial_{x}w^{2}\|_{1,t}^{2}\,dt)^{1/2}\\
\|\int_{0}^{t}U(t-t')\partial_{x}w^{2}(t')dt'\|_{8}&=&
\|\partial_{x}\int_{0}^{t}P_{-}U(t-t')(D_{x}^{\sigma_{1}}
D_{y}^{\gamma_{1}})\partial_{x}(w^{2})(t')
dt'\|_{L^{\infty}_{y}L^{2}_{x,t}}\\
&\lesssim&\|\partial_{x}(D_{x}^{\sigma_{1}}
D_{y}^{\gamma_{1}})(w^{2})\|_{L^{1}_{y}L^{2}_{x,t}}
\lesssim (\int\|\partial_{x}w^{2}\|_{2,t}^{2}dt)^{1/2}\\
\|\int_{0}^{t}U(t-t')\partial_{x}w^{2}(t')dt'\|_{9}&\lesssim&
\int_{0}^{T}\|\partial_{x}w^{2}(t')\|_{2,t'}dt'+
\int_{0}^{T}\|\partial_{x}w^{2}(t')\|_{1,t'}dt'\\
\|\int_{0}^{t}U(t-t')\partial_{x}v^{2}(t')dt'\|_{10}&\lesssim&
\int_{0}^{T}\|\partial_{x}w^{2}(t')\|_{3,t'}dt'+
\int_{0}^{T}\|\partial_{x}w^{2}(t')\|_{1,t'}dt'\\
\|\int_{0}^{t}U(t-t')\partial_{x}v^{2}(t')dt'\|_{11}&\lesssim&
\int_{0}^{T}\|\partial_{x}w^{2}(t')\|_{3,t'}dt'+
\int_{0}^{T}\|\partial_{x}w^{2}(t')\|_{4,t'}dt'+
\int_{0}^{T}\|\partial_{x}w^{2}(t')\|_{1,t'}dt'\\
\|\int_{0}^{t}U(t-t')\partial_{x}v^{2}(t')dt'\|_{13}&\lesssim&
\int_{0}^{T}\|\partial_{x}w^{2}(t')\|_{4,t'}dt'+
\int_{0}^{T}\|\partial_{x}w^{2}(t')\|_{3,t'}dt'+
\int_{0}^{T}\|\partial_{x}w^{2}(t')\|_{1,t'}dt'.
\end{eqnarray*}
It is now clear that the heart of the matter is reduced
to obtaining  good estimates for
$$
\|\langle y\rangle^{\alpha}D_{x}^{2}
(\partial_{x}w^{2})\|_{L^{2}_{x,y,t}},\, \, \,
\|\langle y\rangle^{\alpha}(D_{y}^{2})
\partial_{x}w^{2}\|_{L^{2}_{x,y,t}},\, \, \,
\|\langle y\rangle^{\alpha}D_{x}^{\sigma}
D_{y}^{\gamma}\partial_{x}w^{2}\|_{L^{2}_{x,y,t}}.
$$
\vskip .2 in
\noindent
{\bf{ Estimate of $\|\langle y\rangle^{\alpha}D_{x}^{2}
(\partial_{x}w^{2})\|_{L^{2}_{x,y,t}}$}}\\
we can write
\begin{eqnarray*}
\|\langle y\rangle^{\alpha}D_{x}^{2}(\partial_{x}w^{2})\|
_{L^{2}_{x,y,t}}&\sim&  \|
\langle y\rangle^{\alpha}\partial_{x}^{3}(w)w\|
_{L^{2}_{x,y,t}}+ \|
\langle y\rangle^{\alpha}\partial_{x}^{2}(w)\partial_{x}w\|
_{L^{2}_{x,y,t}} \\
&\sim&\|
\langle y\rangle^{\alpha}\partial_{x}^{3}(P_{+}w)w\|
_{L^{2}_{x,y,t}}+
\|\langle y\rangle^{\alpha}\partial_{x}^{3}(P_{-}w)w)\|
_{L^{2}_{x,y,t}}+\|\langle y\rangle^{\alpha}\partial_{x}^{2}(w)
\partial_{x}w\|_{L^{2}_{x,y,t}} \\
&\lesssim& \|w\|_{5}\|w\|_{10} +\|w\|_{7}\|w\|_{9}+
\|w\|_{1}\|w\|_{12}.
\end{eqnarray*}

\vskip .2 in
\noindent
{\bf{ Estimate of $\|\langle y\rangle^{\alpha}(D_{y}^{2})
\partial_{x}w^{2}\|_{L^{2}_{x,y,t}}$}}\\
With similar arguments and (\ref{weight}) one obtains
$$
\|\langle y\rangle^{\alpha}(D_{y}^{2})
\partial_{x}w^{2}\|_{L^{2}_{x,y,t}}\lesssim
\|w\|_{5}\|w\|_{10} +\|w\|_{7}\|w\|_{9}+\|w\|_{1}(\|w\|_{12}
+\|v\|_{11}).
$$

\vskip .2 in
\noindent
{\bf Estimate of $\|\langle y\rangle^{\alpha}D_{x}^{\sigma}
D_{y}^{\gamma}\partial_{x}w^{2}\|_{L^{2}_{x,y,t}}$}\\
Assume now that $\sigma+\gamma\leq 2$. Then
from (\ref{ygroup}) and (\ref{weight}) it
follows that
\begin{eqnarray*}
\|\langle y\rangle^{\alpha}D_{x}^{\sigma}
D_{y}^{\gamma}\partial_{x}w^{2}\|_{L^{2}_{x,y,t}}
   &\lesssim&\|D_{x}^{\sigma}
   D_{y}^{\gamma}
   (\langle y\rangle^{\alpha}\partial_{x}w^{2})\|_{L^{2}}\\
&\lesssim&\|D_{x}^{2}
   (\langle y\rangle^{\alpha}\partial_{x}w^{2})\|_{L^{2}}
   +\|D_{y}^{2}
   (\langle y\rangle^{\alpha}\partial_{x}w^{2})\|_{L^{2}}\\
&\lesssim&\|\langle y\rangle^{\alpha}D_{x}^{2}
   (\partial_{x}w^{2})\|_{L^{2}}
   +\|\langle y\rangle^{\alpha}D_{y}^{2}
   (\partial_{x}w^{2})\|_{L^{2}}
\end{eqnarray*}
and here we use the previous two steps.

We return to the estimates for the
terms in \eqref{l1} involving $v$. It is not hard to show that
estimates similar to those presented above are available for the term
with $wv$ as long as $v\in Y_{T}$. Estimating the term involving
$v^{2}$ is even easier because
$$QU(t-t')\partial_{x}(v^{2})\sim U(t-t')(v^{2}),$$
that is the derivative is ``inactive'' in this part of the
operator $L_{1}$. We then obtain that
\begin{equation}
\|L_{1}(w,v)\|_{X}\lesssim  \| w_{0}\|_{X,0}+
\|w\|_{X}^{2} +\|v\|_{Y}^{2}.
\label{fp1}\end{equation}
Notice that no factor of $T$ appears in front of
$(\|w\|_{X}^{2} +\|v\|_{Y}^{2})$. This
will force us to assume that $\| w_{0}\|_{X,0}$
is small in order to claim that the operator $L$
is indeed a contraction. We will return to this matter a bit later.
We start the estimate for the operator $L_{2}(w,v)$.
By unitarity of the linear flow in $L^2$ Sobolev spaces, we have
$$\|\chi_{[0,T]}U(t)(Id -Q)v_{0}\|_{1}\lesssim \|v_{0}\|_{H^{2}}.$$
Next, from (\ref{ygrouplow})
$$\|\chi_{[0,T]}U(t)( Id -Q)v_{0}\|_{3}\lesssim
\|\langle y\rangle^{\alpha}D_{y}^{\gamma_{2}}v_{0}\|_{L^{2}}
+\|D_{x}^{-\alpha-}D_{y}^{\gamma_{2}+\alpha}v_{0}\|_{L^{2}}.$$
From (\ref{ytmaxxlow}) we have
$$\|\chi_{[0,T]}U(t)( Id -Q)v_{0}\|_{9}\lesssim
\|\langle y\rangle^{\alpha}D_{x}^{-\gamma}
(1+D_{y})^{\gamma}v_{0}\|_{L^{2}}
+\|D_{x}^{-\sigma}(1+D_{y})^{\theta}v_{0}\|_{L^{2}}$$
where  $\gamma_{1}=1/2+, \sigma=1/4+, \theta=1+$.
From (\ref{ytmaxylow}), we have
\begin{eqnarray*}
   \|\chi_{[0,T]}U(t)(Id -Q)v_{0}\|_{10}&\lesssim&
\|\langle y\rangle^{\alpha}
D_{x}^{-\sigma}(1+D_{y})^{\theta}v_{0}\|_{L^{2}}
+\|D_{x}^{-\sigma-1/2}(1+D_{y})^{\theta+1/2}v_{0}\|_{L^{2}},
\end{eqnarray*}
From (\ref{yxyinfinitylow}), we also have
$$\|\chi_{[0,T]}U(t)(Id -Q)v_{0}\|_{11}\lesssim
\|\langle y\rangle^{\alpha}(1+D_{y})^{1+\epsilon}v_{0}\|_{L^{2}}
+\|D_{x}^{-\alpha-\epsilon}(1+D_{y})^{1+\alpha+2\epsilon}
v_{0}\|_{L^{2}},$$
and we ask that $1+\epsilon\leq \gamma_{2}$. Similarly,
$$\|\chi_{[0,T]}U(t)(Id -Q)v_{0}\|_{12}\lesssim
\|\langle y\rangle^{\alpha}(1+D_{y})^{\epsilon}v_{0}\|_{L^{2}}
+\|(1+D_{y})^{\alpha+2\epsilon}v_{0}\|_{L^{2}}.$$
Now we pass to the estimate of the nonlinear terms. For any
$i=1,3,9,10,11,12$ and smooth $f$ and $h$ we have
\begin{eqnarray*}
   \|\chi_{[0,T]}\int_{0}^{t}
U(t-t')(\partial_{x}(fh)dt'\|_{i,t'}&\lesssim&
\int_{0}^{T}\|U(t-t')(\partial_{x}(fh)(t')\|_{i}dt'\\
&\lesssim&T^{1/2}\|(D^{2}_{x}+D_{y}^{2})(fh)(t')\|_{L^{2}_{t,x,y}}
+T^{1/2}\|\langle y\rangle^{\alpha}D_{y}^{\gamma_{2}}
(fh)(t')\|_{L^{2}_{t,x,y}}\\
&\lesssim&T^{1/2}(\|f\|_{1}\|h\|_{11}).
\end{eqnarray*}
This and  Definition \ref{spacexy} show that
\begin{equation}
\|L_{2}(w,v)\|_{X}\lesssim  \| v_{0}\|_{Y,0}+
T^{1/2}(\|w\|_{X}^{2} +\|v\|_{Y}^{2}).
\label{fp2}\end{equation}
Then, by \eqref{fp1} and \eqref{fp2}, it is easy to conclude that
$L$ is a contraction (and hence has a unique fixed point),
if we assume \eqref{smalldelta},
for some small $\delta$. This does not
yet prove Theorem \ref{main1} and \ref{main2}. But before we present
the  conclusion of the proofs we would like to make few remarks.
\begin{remark}
Based on the order of derivatives needed
for the maximal function estimates in Theorem \ref{max} and
\ref{ymax}, one may  guess that a local well-posedness
theory could be obtained in a Sobolev space of order at most
$3/2+\epsilon$. In fact,
one can prove that $\|[D_{y}^{3/2+}\partial_{x}(u)]
[\langle y\rangle^{\alpha}u]\|_{L^{2}_{x,y,t}}$ can be bounded using
the smoothing effect norms of type
$\|\cdot\|_{i},$ for $ i=5,6,7,8,$ of  order at most $3/2+$,
and maximal function norms of type
$\|\cdot\|_{j},$ for $ j=9,10$ . But, unfortunately, the intermediate
terms appearing in the Leibniz rule cannot be handled in a simple way.
Take for example the expression $\|D_{y}^{\gamma}(D_{y}\partial_{x}(u)
u\langle y\rangle^{\alpha})\|_{L^{2}_{x,y,t}}, \gamma=1/2+$.
To estimate this term one would need
interpolations norms between $\|\cdot\|_{i},$ for $ i=5,6,7,8,$ of
order at most $3/2+\epsilon$ and
$\|\cdot\|_{j},$ for $ j=9,10$. Unfortunately these are not
available when
the two functions involved are $D_{y}\partial_{x}P_{\pm}(u)$
and $P_{\mp}u$. One would need to prove directly
the estimates for the intermediate norms of the linear solution
$U(t)u_{0}$. We decided not to attempt this here because we believe
that the well-posedness results that one can obtain from the
oscillatory integral theory would not be  optimal
(compare for example with \cite{KPV1} and \cite{KPV4}).
\end{remark}
\begin{remark}
In the estimate of $L_{2}$, the low-frequency part of the
operator $L$, a power $T^{1/2}$ appears.
This is because we do not need  the
``smoothing effect'' norm, the only norm that prevents a factor
$T^{\epsilon}$ from appearing. This also says that if one is
interested only on a local well-posedness result, the smallness
assumption on $\|v_{0}\|_{Y,0}$ can be removed.
\label{lwp}\end{remark}
Now we go back to the conclusion of the proof for Theorem
\ref{main1} and \ref{main2}. We need to relax the
 smallness assumption \eqref{smalldelta}.
We will prove that for global well-posedness we only need
the norms
$$\|\langle y\rangle^{\alpha}w_{0}\|_{L^{2}}, \, \,
\|\langle y\rangle^{\alpha}D_{x}^{-1/2-\epsilon}v_{0}\|_{L^{2}}, \, \,
\|D_{x}^{-3/4-\epsilon}v_{0}\|_{L^{2}}$$
to be small. To do so, we
rescale the solution $u$ by observing that if $u(x,y,t)$ solves
(\ref{ivp}) on $[0,T]$, then $u_{\rho}(x,y,t)=
\rho^{2}u(\rho x,\rho^{2}y,\rho^{3}t)$ also solves (\ref{ivp})
with initial data $u_{0,\rho}(x,y)=\rho^{2}u_{0}(\rho x,\rho^{2}y)$,
on the interval $[0,\rho^{-3}T]$.
We have the following lemma:
\begin{lemma}
Assume $\rho \in [0,1)$ and
$v_{\rho}(x,y)=\rho^{2}u(\rho x,\rho^{2}y)$. Then for any
$\sigma, \gamma\geq 0$
\begin{eqnarray}
    \label{riscaling1}
\|D_{x}^{\sigma}u_{\rho}\|_{L^{2}}&\lesssim&\rho^{1/2+\sigma}
\|D_{x}^{\sigma}u\|_{L^{2}}\\
   \label{riscaling2} \|D_{y}^{\gamma}u_{\rho}\|_{L^{2}}&\lesssim&
\rho^{1/2+2\gamma}\|D_{y}^{\gamma}u\|_{L^{2}}\\
\label{yriscaling1} \||y|^{\alpha}D_{x}^{\sigma}
u_{\rho}\|_{L^{2}} &\lesssim&
    \rho^{1/2+\sigma-2\alpha}
    \||y|^{\alpha}D_{x}^{\sigma}u\|_{L^{2}}.
\end{eqnarray}
If then $\gamma>1$, we also have
\begin{equation}
\label{yriscaling2} \||y|^{\alpha}D_{y}^{\gamma}
u_{\rho}\|_{L^{2}} \lesssim
    \rho^{1/2+2\gamma-2\alpha}\|<y>^{\alpha}D_{y}^{\gamma}u
    \|_{L^{2}}.
\end{equation}
\label{risc}\end{lemma}
\begin{proof}
First we observe that
$$
\widehat{D_{y}^{\gamma}u_{\rho}}(\x,\lambda)=
\rho^{-1}|\lambda|^{\gamma}\hat{u}(\rho^{-1}\x,\rho^{-2}\lambda).
$$
We then have
\begin{eqnarray*}
 \|D_{y}^{\gamma}u_{\rho}\|_{L^{2}}^{2}&=&
 \||\lambda|^{\gamma}\widehat{u_{\rho}}\|_{L^{2}}^{2}\\
 &=&\rho^{-2}\int |\lambda|^{2\gamma}|\hat{u}|^{2}
 (\rho^{-1}\x,\rho^{-2}\lambda)\,d\x\,d\lambda\\
 &\lesssim&\rho^{1+2\gamma}\int |\lambda|^{2\sigma}|\hat{u}|^{2}
 (\x,\lambda)\,d\x\,d\lambda \lesssim \rho^{1+2\gamma}
 \|D_{y}^{\gamma}u\|_{L^{2}}^{2}.
\end{eqnarray*}
and (\ref{riscaling2}) follows. The same argument can be used to
prove (\ref{riscaling1})
We now prove (\ref{yriscaling1})
by complex interpolation. If $z=i\beta$, then by  (\ref{riscaling1})
$$\||y|^{z}D_{x}^{\sigma}u_{\rho}\|_{L^{2}}
\lesssim \rho^{1/2+2\sigma}\|D_{x}^{\sigma}u\|_{L^{2}}.$$
If $z=1+i\beta$ then by (\ref{riscaling1})
\begin{eqnarray*}
\||y|^{z}D_{x}^{\sigma}u_{\rho}\|_{L^{2}}&=&
\||\xi|^{\gamma}\partial_{\lambda}\widehat{u_{\rho}}
\|_{L^{2}}\\
&=&\rho^{-2}\|D_{x}^{\sigma}(yu)_{\rho}\|_{L^{2}}\\
&\lesssim& \rho^{1/2+2\sigma-2}(\|yD_{x}^{\sigma}u\|_{L^{2}}
\end{eqnarray*}
Interpolation gives immediately (\ref{yriscaling1}).
We are now ready to prove (\ref{yriscaling2}). Again we use
interpolation. As above, for $z=i\beta$,
$$\||y|^{z}D_{y}^{\gamma}u_{\rho}\|_{L^{2}}
\lesssim \rho^{1/2+2\gamma}\|D_{y}^{\gamma}u\|_{L^{2}}.$$
If $z=1+i\beta$ then
\begin{eqnarray*}
\||y|^{z}D_{y}^{\gamma}u_{\rho}\|_{L^{2}}&=&
\|\partial_{\lambda}(|\lambda|^{\gamma}\widehat{u_{\rho}})
\|_{L^{2}}\\
&=&\||\lambda|^{\gamma-1}\widehat{u_{\rho}}\|_{L^{2}}+
\||\lambda|^{\gamma}\partial_{\lambda}(\widehat{u_{\rho}})
\|_{L^{2}}\\
&=&\|D_{y}^{\gamma-1}u_{\rho}\|_{L^{2}}+
\rho^{-2}\|D_{y}^{\gamma}(yu)_{\rho}\|_{L^{2}}\\
&\lesssim& \rho^{1/2+2\gamma-2}(\|D_{y}^{\gamma-1}u\|_{L^{2}}+
\|D_{y}^{\gamma}(yu)\|_{L^{2}}).
\end{eqnarray*}
To finish the chain of inequalities we use \eqref{weight} and we
obtain
$$\||y|^{z}D_{y}^{\gamma}u_{\rho}\|_{L^{2}}\lesssim
(\rho^{1/2+2\gamma-2})
\|\langle y\rangle D_{y}^{\gamma}u\|_{L^{2}}.$$
Interpolation gives immediately (\ref{yriscaling2}).

\end{proof}
It is now easy to see from  (\ref{riscaling1})-(\ref{yriscaling2})
that if we repeat the fixed point argument above for\\ $Qu_{\rho}=
u_{1,\rho}$ and $(Id-Q)u_{\rho}=
u_{2,\rho}$, and
$$\max(\|\langle y\rangle^{\alpha}w_{0,\rho}\|_{L^{2}}, \, \,
\|\langle y\rangle^{\alpha}D_{x}^{-1/2-\epsilon}v_{0,\rho}\|_{L^{2}}, \, \,
\|D_{x}^{-3/4-\epsilon}v_{0,\rho}\|_{L^{2}})\ll \delta$$
then $\max(\|w_{0,\rho}\|_{X,0}, \|v_{0,\rho}\|_{X,0})
\leq \delta$ and the fixed point
argument can be applied for an appropriate small $\delta$. This
concludes the proof of Theorem~\ref{main1}. To conclude also the
proof of Theorem~\ref{main2} one needs to use Remark \ref{lwp}
and
(\ref{riscaling1})-(\ref{yriscaling2}), as we did above.
\begin{remark}
One can also consider the modified KP-I initial value problem
\begin{equation}
\left\{ \begin{array}{l}
\partial_x(\partial_t u + \partial_x^3 u
+\beta u^{2}\partial_xu)+ \partial_y^2u= 0, \\
u(x,0) = u_0(x), \hspace{1.5cm}(x,y) \in \dbR^2, \,  t \in \dbR.
\end{array}\right.
\label{ivp0}\end{equation}
Using the arguments presented above, one can prove a local
well-posedness theorem in $H^{2}$, without introducing spaces
involving the weight $\langle y\rangle^{\alpha}$ and without assuming
smallness of the initial data. In fact the weight
$\langle y\rangle^{\alpha}$
was introduced above to transform a norm $L^{1}_{y}$, needed in the 
inhomogeneous smoothing effect (\ref{nhp-}),  into an $L^{2}$
norm needed for the maximal function estimates (\ref{tmaxx}) and
(\ref{tmaxy}). In modified KP-I 
this is not needed because
the square power in the nonlinearity $u^{2}\partial_{x}u$
takes us for free from the $L^{1}$ norm to the $L^{2}$.
Then, based on (\ref{riscaling1}) and  (\ref{riscaling2})
we can rescale the solution and remove the smallness assumption.
The precise well-posedness statement for  \eqref{ivp0} can be
summarized in the theorem below. Define the space $W_{0}$ through the
norm
$$\|f\|_{W_{0}}=\|f\|_{H^{2}}+\|D_{x}^{-\gamma}
(1+D_{y})^{\gamma}f\|_{L^{2}}+\|D_{x}^{-\sigma}
(1+D_{y})^{\theta}f\|_{L^{2}},$$
where $\gamma>1/2, \sigma>1/4, \theta>1$.
\begin{theorem}
For any $u_{0}\in W_{0}$ there exist $T=T( \|u_{0}\|_{W_{0}})$
and a unique solution $u$ for \eqref{ivp0} such that
$u(x,y,t) \in C([0,T],W_{0})$. Moreover the map
$u_{0}\longrightarrow u$ is continuous with
respect to the initial  data in the appropriate topology.
\end{theorem}
An adaptation of the arguments from \cite{KPV2} to the KP-I setting
should, in principle, provide similar local well-posedness results
for pure power generalizations (with nonlinearity $u^p u_x,~p \in 
{\mathbb{N}}$) of KP-I. These extensions would lower the regularity
required for existence in the energy method argument of Tom \cite{T}
and also provide uniqueness. 
\end{remark}

\section{Appendix: The fractional  Leibniz Rule}
In this section we recall some known facts on fractional
Leibniz rule for one variable functions and some related results
involving the weight $\langle y\rangle^{\alpha}$.

\begin{theorem}
Assume $0<\sigma<1$ and $1<p<\infty$. Then
$$\|D^{\sigma}_{x}(fg)-fD^{\sigma}_{x}(g)
-gD^{\sigma}_{x}(f)\|_{L^{p}_{x}}
\lesssim \|g\|_{L^{\infty}}\|D^{\sigma}_{x}(f)\|
_{L^{p}}.$$
\label{linfinity}\end{theorem}
For the proof one can see \cite{KPV2}.
We  also need  a lemma that relates fractional
derivatives with  the weight $\langle y\rangle^{\alpha}$.
\begin{lemma}
Assume $0\leq \alpha\leq 1, 0<\gamma<1$ and $1<p<\infty$. Then
\begin{equation}
   \|\langle y\rangle^{\alpha}(1+D)^{\gamma}_{y}(f)
   \|_{L^{p}_{y}}=
   \|(1+D)^{\gamma}_{y}(\langle y\rangle^{\alpha}f)
   \|_{L^{p}_{y}}+
   O(\|\langle y\rangle^{\alpha}f\|_{L^{p}_{y}}).
\label{weight}\end{equation}
\end{lemma}
\begin{proof}
It is enough to prove the following estimate on the commutator:
\begin{equation}
\|[\langle y\rangle^{\alpha},(1+D_{y})^{\gamma}]f\|_{L^{p}}\lesssim
\|\langle y\rangle^{\alpha}f\|_{L^{p}_{y}}.
\label{part1}\end{equation}
we start by first proving that for any $\beta \in \dbR$
\begin{equation}
\|D^{\gamma}(\langle y\rangle^{i\beta})\|
_{L^{\infty}}\leq C.
\label{hope}\end{equation}
By Proposition 1 page 241 in \cite{Stein} we have that
\begin{equation}
   |K_{\beta}(\lambda)|=|\F(\langle y\rangle^{i\beta})(\lambda)|
   \lesssim |\lambda|^{-1-N},
\label{kernel}\end{equation}
for all $N\geq 0$. Then
\begin{eqnarray*}
   |D^{\gamma}(\langle y\rangle^{i\beta})|&=&
\left|\int e^{-iy\lambda}|\lambda|^{\gamma}
K_{\beta}(\lambda)\, d\lambda\right|\\
&\lesssim& \int_{|\lambda|<1}|\lambda|^{\gamma-1}\, d\lambda+
\int_{|\lambda|\geq 1}|\lambda|^{\gamma-3}\, d\lambda<\infty,
\end{eqnarray*}
where in the last step we used (\ref{kernel}) with $N=0$ and $N=2$.
We are now ready to prove (\ref{part1}) by complex
interpolation. For $z \in \dbC, 0\leq\Re\,z\leq 1$ we define
the operator $T^{z}=[\langle y\rangle^{2z},(1+D)^{\gamma}]$.
Assume now $z=i\beta$. Then by Theorem \ref{linfinity}
\begin{eqnarray*}
\|T^{z}f\|_{L^{2}}&=&\| \langle y\rangle^{i2\beta}(1+D)^{\gamma}f
-(1+D)^{\gamma}(\langle y\rangle^{i2\beta}f)\|_{L^{p}}\\
&\lesssim&\|[(1+D)^{\gamma}-D^{\gamma}]f\|_{L^{p}}
+\|[(1+D)^{\gamma}-D^{\gamma}]\langle y\rangle^{i2\beta}f\|_{L^{p}}\\
&+&\| \langle y\rangle^{i2\beta}D^{\gamma}f
-D^{\gamma}(\langle y\rangle^{i2\beta}f)+
D^{\gamma}(\langle y\rangle^{i2\beta})f\|_{L^{p}}
+\|D^{\gamma}(\langle y\rangle^{i2\beta})f\|_{L^{p}}\\
&\lesssim&\|[(1+D)^{\gamma}-D^{\gamma}]f\|_{L^{p}}
+\|[(1+D)^{\gamma}-D^{\gamma}]\langle y\rangle^{i2\beta}f\|_{L^{p}}\\
&+&\|(1+D)^{\gamma}(\langle y\rangle^{i2\beta})\|
_{L^{\infty}}\|f\|_{L^{p}}+\|f\|_{L^{p}}\lesssim \|f\|_{L^{p}}
\end{eqnarray*}
where in the last step we used (\ref{hope}) and the fact that
the multiplier $m(\lambda)=
(|\lambda|+1)^{\gamma}-|\lambda|^{\gamma}$ is a good Marcinkiewicz
multiplier (see \cite{Stein} page 245).
Assume now that $z=1+i\beta$. It is enough to estimate
$\tilde{T}^{z}f=\langle y\rangle^{2i\beta}y^{2}(1+D)^{\gamma}f
-(1+D)^{\gamma}(y^{2}\langle y\rangle^{i2\beta}f)$.
Observe that
$$y^{2}(1+D^{\gamma})f=M_{2}f+M_{1}yf+(1+D)^{\gamma}(y^{2}f)$$
where $M_{j}$ are multiplier operators such that
$$\F(M_{j}f)(\lambda)=(1+|\lambda|)^{\gamma-j}
\hat{f}(\lambda).
$$
Using this representation, Theorem \ref{linfinity}, (\ref{hope})
and and again  theory of multiplier operators, we have
\begin{equation}
\|\tilde{T}^{1+i\beta}f\|_{L^{p}}\sim
\|f\|_{L^{p}}+ \|yf\|_{L^{p}}+ \|y^{2}f\|_{L^{p}}\lesssim
\|(1+y^{2})f\|_{L^{p}},
\end{equation}
Then complex interpolation gives (\ref{part1}).
\end{proof}

\noindent
{\bf Acknowledgment.} G. Staffilani thanks the Department of
Mathematics at the University of Chicago for the generous
hospitality  during the preparation of this paper.

\end{document}